\input epsf
\title{
Reconstruction of the intertwining operator\\ 
and new striking examples added to\\ 
``Isospectral pairs of 
metrics on balls and spheres 
with different local geometries"
}

\newcommand{\cmt}[1]{\ifhmode\newline\fi{\sf *** \ \ #1 \\}}

\documentclass[11pt]{article}
\usepackage{latexsym}
\usepackage{epsfig}

\usepackage{amsmath}

\usepackage{amssymb}

\usepackage{exscale}

\usepackage{amsthm}
\usepackage{verbatim}

\newtheorem{theorem}{Theorem}[section]
\newtheorem{lemma}[theorem]{Lemma}

\newcommand{\R}{{\mathbf{R}}}

\def\:{\colon}

\long\def\onefigure#1#2{
\begin{figure*}[tbh]
\begin{center}
#1
\end{center}
\caption{#2}
\end{figure*}
} 

\newcommand{\iipefig}[1]  
{\smallskip\begin{center}\def\IPEfile{#1.ipe}\input{\IPEfile}\end{center}\smallskip}
\newcommand{\lipefig}[2]  
{\onefigure{\def\IPEfile{frh-#1.ipe}\input{\IPEfile}}{\label{f:#1} #2} }

\author
{Z. I. Szab\'o\thanks{Partially supported by
NSF grant DMS-0604861}
\thanks{Lehman College of CUNY, Bronx, NY 10468,USA, and 
R\'enyi Institute, Budapest,
POBox 127, 1364 Hungary. E-mail: zoltan.szabo@lehman.cuny.edu}}
\date{}
\begin{document}

\maketitle


\begin{abstract}
The intertwining operator constructed in
\cite{sz1,sz2}  
\footnote{Annals of Mathematics, 154(2001), 437-475; 
and 161(2005), 343-395}
does not appear in the right form. It is established there by
using only the 
anticommutators $J_1$ and $J^\prime_1$. The correct operator   
involves all endomorphisms, $J_\alpha$, which are unified
by the Z-Fourier transform. 
Although some of the correct elements of the
previous constructions are kept, this idea is
established by a new technique 
which yields the various
isospectrality theorems stated in the papers 
on a much larger scale.
The new results include   
new isospectrality examples living on sphere$\times$ball-
and sphere$\times$sphere-type manifolds. Among them, there are such
discrete isospectrality families where one of the members is
homogeneous while the others are locally inhomogeneous (striking examples).
Furthermore, a large class of new isospectrality families are constructed
by $\sigma$ deformations.
\end{abstract}

\section{Introduction.} 
In papers \cite{sz1,sz2}, 
the intertwining operator is constructed 
by the complex linear correspondence 
\begin{eqnarray}
\label{illkappa}
\kappa^* :
\varphi (|X|,Z)\Theta_{Q_1}(X,J_1)...\Theta_{Q_p}(X,J_1)\overline
\Theta_{Q_{p+1}}(X,J_1)\dots\overline\Theta_{Q_{p+q}}(X,J_1)
\\
\to
\varphi (|X|,Z)\Theta_{Q_1}(X,J^\prime_1)\dots
\Theta_{Q_p}(X,J^\prime_1)\dots\overline
\Theta_{Q_{p+1}}(X,J^\prime_1)\dots\overline\Theta_{Q_{p+q}}(X,J^\prime_1),
\nonumber
\end{eqnarray}
where $Q_1,\dots ,Q_{p+q}$ are arbitrary X-vectors, furthermore, 
$\Theta_Q(X,J_1)=\langle Q+\mathbf iJ_1(Q),X\rangle$ and
the corresponding $\Theta^\prime_Q(X,J^\prime_1)$ are defined
by the anticommutators 
$J_1$ and $J^\prime_1$ respectively. 
These anticommutators
$\sigma$-relate to each other, meaning, that the X-space is
a direct sum, 
$
\mathbf v=
\mathbf v^{(a)}
\oplus
\mathbf v^{(b)}
$,
of 
$
k^{(a)}
$
resp.
$
k^{(b)}
$
real-dimensional subspaces
such that the components are invariant under the action 
of the anticommutators and they agree on the first component while 
$J_1=-J^\prime_1$ holds on the second one. 
The complex linear property
guaranties that it intertwines the operators 
$\partial_1D_1\bullet$ and $\partial_1D^\prime_1\bullet$,  
which appear in the Laplacians involved to the angular momentum operators 
$\mathbf M=\sum_{\alpha =1}^l\partial_\alpha D_\alpha\bullet$ resp. 
$\mathbf M^\prime$. 
One of the ultimate goals is to intertwine the complete Laplacians.

It was overlooked by the author that one can not allow arbitrary 
vectors, $Q_i$, in the
definition of  
$\kappa^*$
because the operator becomes ill-defined. In fact, well-defined 
complex linear map can be introduced by choosing a system,
$\{E_i^{(a)} ,E^{(b)}_{j}\}$, where 
$1\leq i\leq
k^{(a)}/2
$;
$1\leq j\leq
k^{(b)}/2
$,
of independent vectors in the corresponding component spaces
$
\mathbf v^{(a)}
$
and 
$
\mathbf v^{(b)}
$
which form a complex 
linear basis with respect
to both complex structures $J_1$ and $J^\prime_1$. Then, the admissible 
$Q$'s are the vectors laying in the real 
subspace spanned by all these $E_r^{(c)}$'s. This natural complex 
linear map can be 
described in terms of the complex coordinates 
$\{z_r\}$ resp. $\{z^\prime_r\}$
determined by the same basis
$
\{E_r|
1\leq r\leq k^{(a)}/2+k^{(b)}/2=k/2
\}=\{E^{(a)}_i\}\cup\{ E^{(b)}_{j}\}
$,
for the two complex structures respectively such
that the image of a polynomial written up in terms of the coordinates
$\{z_r\}$ is the polynomial of the same form but written 
up in terms of the other coordinates $\{z^\prime_r\}$.
One can easily see that these
well-defined maps depend on the real subspaces spanned by  
$\{E_r\}$, meaning that maps defined for different real
subspaces correspond different elements 
to the very same element in general.
Thus, by allowing all possible $Q$'s, the above $\kappa^*$
is ill-defined indeed. 

Unfortunately, this problem can not be eliminated by replacing the 
ill-defined operator by one of these well-defined ones. 
In this case a much more serious 
difficulty appears, namely, such a well-defined operator 
does not intertwine 
$\partial_\alpha D_\alpha\bullet$ and 
$\partial_\alpha D^\prime_\alpha\bullet$ 
satisfying $\alpha >1$, which are also parts 
of the corresponding Laplacians. 
Let it also be mentioned that the rest parts of the 
Laplacians as well as the boundary conditions are intertwined by it.
In the papers, the proof of intertwining of the above parts of 
$\mathbf M$ resp. $\mathbf M^\prime$
explores the false assumption claiming that the $\kappa^*$ operates
on functions defined for imaginary $Q$'s in the same way as for
the real ones.  

One of the reasons causing this blunder was that, instead of 
$\mathbf M$ and $\mathbf M^\prime$, one was focusing just on 
$D_\alpha\bullet$
and $D_\alpha^\prime\bullet$, i. e., tried to define intertwining 
operator using only the X-space. An other reason was that the intertwining 
operator was defined just with endomorphisms 
$J_1$ and $J^\prime_1$.
The corrected operator, involving all $J_\alpha$ and focusing on 
$\mathbf M$ and $\mathbf M^\prime$, is defined by choosing a basis 
$\{E_1,\dots ,E_{k/2}\}$ and using $Q_i$'s laying in the real span of these
basis vectors. More precisely, this correspondence is:
\begin{eqnarray}
\label{goodkappa}
\kappa :
\int_{\mathbf z} e^{\mathbf i\langle Z,V\rangle}
\varphi (|X|,V)
\Theta_{Q_{1\dots p}}
(X,V_u)
\overline\Theta_{Q_{p+1\dots p+q}}
(X,V_u)
dV
\to
\\
\int_{\mathbf z} e^{\mathbf i\langle Z,V\rangle}\varphi (|X|,V)
\Theta^\prime_{Q_{1\dots p}}
(X,V_u)
\overline\Theta^\prime_{Q_{p+1\dots p+q}}(X,V_u)
dV,
\nonumber
\end{eqnarray}
where
$
\Theta_{Q_{1\dots p}}(X,V_u):=
\Theta_{Q_1}(X,V_u)\dots
\Theta_{Q_p}(X,V_u)
$
and the corresponding
$
\Theta^\prime_{Q_{1\dots p}}(X,V_u)
$
are defined by the endomorphisms $J_{V_u}$ and $J^\prime_{V_u}$ 
belonging to the unit Z-vectors $V_u$ respectively.
This operator associates functions defined by the Z-Fourier transform 
to each other. Appropriate intertwining, $\kappa^-$, can be established
by using $e^{-\mathbf i\langle Z,V\rangle}$ in the above formulas. 
Then also 
$\kappa^{\mathbf R}=(\kappa +\overline{\kappa}^-)/2$ is going to be an 
intertwining operator. However, this paper proceeds only
with the first version. 

The well-definedness of this operator follows from reasons such as only
$Q$'s laying in the real span of vectors $E_i$ are used in its definition,
furthermore, the Fourier transform is an isometry on the corresponding
$L^2$ Hilbert spaces. The intertwining of the
Laplacians is due to the Z-Fourier transform implemented into the formulas. 
The addition of this Z-Fourier transform to the original idea makes 
the mathematical situation much more complex, requiring a complete 
rethinking of the original construction. For instance,
beyond proving the intertwining
of the Laplacians, it is much more difficult to prove the intertwining
of the boundary conditions for this operator. Actually, the proof of the 
latter statement combines the Z-Fourier transform with an independent 
idea incorporated into the {\it Independence Theorems}. 
The ill-definedness of 
(\ref{illkappa})
was recognized by H. F\"urstenau whose observation triggered
the author's thorough rethinking of his complete construction. 
The much deeper problems hidden 
under the cover of ill-definedness came to the light during this revision
process. The reborn 
operator presented here saves all the previous results and 
provides also new interesting isospectrality examples.
The main goal is to establish  

\begin{theorem}
\label{th_intr}
The ball- resp. sphere-type manifolds, which have the same radius function
$\varphi (|X|,|Z|)$ and are defined   
on H-type groups $H^{(a,b)}_l$ having the same parameters $a+b$ and $l$,
are isospectral. This statement extends to a large class of general 
2-step nilpotent Lie groups where an isospectrality family is defined 
by $\sigma$-deformations. 

On ball-type manifolds this statement includes the well-definedness of
(\ref{goodkappa}) and the intertwining of the Laplacians as well
as the boundary conditions.
Since the Dirichlet condition is intertwined, by restrictions, the
operator induces bijections between the function spaces defined
on the boundaries. By observing that it is enough to use functions
satisfying the Z-Neumann condition on the ambient manifolds, one
can prove the intertwining property also on the boundary manifolds.  

All these statements extend onto the solvable extensions of 2-step
nilpotent Lie groups. Furthermore, new examples, not discussed
in the original papers, are also constructed. They live on 
sphere$\times$ball- and
sphere$\times$sphere-type manifolds. Among them two are particularly
interesting. Namely, the isospectrality family of 
sphere$\times$sphere-type manifolds, constructed both on  
$H^{(a,b)}_3$ and
$SH^{(a,b)}_3$, 
the metric is homogeneous for 
the manifold belonging
to the pair $(a+b,0)$ or $(0,a+b)$, while the metrics satisfying 
$ab\not =0$ are locally inhomogeneous. Also the dimensions of
groups of isometries acting on the members are different. These are
new contributions to the old list of striking examples constructed 
on the sphere$\times$torus-type manifolds of 
$H^{(a,b)}_3$. 
resp. geodesic spheres of
$SH^{(a,b)}_3$. 

These theorems are established, first, on H-type groups and their solvable
extensions. Then, they are extended to those 2-step nilpotent groups and 
their solvable extensions which are defined by endomorphism spaces obtained
by perturbing the endomorphism space of a given H-type group
$H^{(a,b)}_l$. 
\end{theorem} 

The perturbation process mentioned above is as follows.
The endomorphisms, $J_Z$, on H-type groups are defined by endomorphisms, 
$j_Z$, acting
on the irreducible components, by the formula
$J_Z=(j_Z,\dots ,j_Z,-j_Z,\dots ,-j_Z)$. The perturbation primarily 
concerns the endomorphisms $j_Z$, i. e., close to the Cliffordian one, a 
new linear space of endomorphisms is chosen with elements denoted by   
$\tilde j_Z$. 
Then, the groups,
$\tilde H^{(a,b)}_l$ 
and
$\tilde{SH}^{(a,b)}_l$, resulted by such a perturbation
arise from endomorphisms 
$\tilde J_Z$ 
defined by the same formula in terms of
$\tilde j_Z$. 

Such a perturbation results a new 
family, defined by the same $(a+b)$ and $l$, whose members are obviously 
$\sigma$-equivalent. By a theorem proved in the last section, the 
$\kappa$ always intertwines the Laplacians, however, the same statement for
the boundary conditions is not guarantied. This problem is solved by the
independence theorems which state the independence of certain subspaces 
formed by 
functions. This independence is established, originally, for H-type groups
defined by Cliffordian endomorphism spaces and remains true for groups
defined by endomorphism
spaces which are close enough to the Cliffordian ones. The latter spaces 
constitute an open set whose members are called the
{\it small perturbations} of a given Cliffordian endomorphism space.
The complete isospectrality theorems are established only on those
$\sigma$-equivalent groups whose endomorphism spaces are produced by small
perturbations of the Cliffordian ones. It is also important
to understand that the independence theorems alone do not validate the
intertwining property for the boundary conditions. They can guarantee 
it just for intertwining operators defined by $\sigma$-deformations.  

The new features in this mathematical process include, first of all,
the integral formula by which the intertwining operator is defined.
This formula deeply roots in quantum theory \cite{sz4,sz5}, furthermore, it 
can be used also for explicit computations of the eigenfunctions and 
spectra. This rooting in physics is exhibited by the surprising fact
that the Laplacian on the investigated manifolds can be identified
with the Landau-Zeeman operator attached to electron-positron systems 
where these particles are orbiting in constant magnetic fields. 
Above, the endomorphisms 
$-j_Z$ resp. $j_Z$
correspond to electrons resp. positrons, and $\sigma$-deformations is
interpreted such
that some of the electrons are exchanged for positrons. Yet, the spectra
on all submanifolds investigated in this paper are not changing during
this exchange-process. The local geometry, however, is dramatically
changing. The manifolds in the striking examples, for instance, are
homogeneous for systems having particles of the same type, while, they are
locally inhomogeneous for mixed particles. 

This is a physical interpretation of the above isospectrality theorems.
The perturbation can be interpreted such that, instead of a system of
identically charged particles, one considers ones which 
are charged distinctly.
Then, the isospectrality theorems are established also for systems produced
by small perturbation of the charge.
The fact of non-changing spectra during electron-positron-exchanges is
well known in physics. However, the statement of this form concerns 
the spectra considered on a non-compact manifold. 
Our statement claims much more than just this. Namely, the spectra remain 
the same also on a large class of compact 
submanifolds. An other distinguishing feature is that no attached local 
geometries are considered in quantum theory.

The intertwining of the Laplacians can be established by using only the
integral formula. For proving the intertwining property for the boundary 
conditions, a new idea, appearing in the Independence Theorems, is 
involved. These theorems are also important new features in this field.

The methods developed in this paper apply only for $\sigma$-deformations.
A characteristic feature of these discrete deformations is that they do
not change the Ricci curvature. This fact is strongly used in proving the
intertwining regarding the boundary conditions. This experience strongly 
indicates that
the submanifolds considered in this paper are not isospectral on
the Gordon-Wilson \cite{gw} examples where continuous isospectral 
deformations with changing Ricci tensor are established on groups 
defined by 2-dimensional Z-spaces. 
Other arguments supporting this statement 
are explained in the end of the last section. 

The much more general Theorem \ref{th_intr} replaces the Isospectrality 
Theorems of the articles. Only the construction of the 
intertwining operator (cf. pages 461-465 in \cite{sz1} and 
371-375 in \cite{sz2}) is effected by this problem. The major 
non-effected part
includes all the Non-Isometry resp. Rigidity Theorems and the
preparatory part of Sections $4.$ resp. $3.$. This
problem with solution was announced at the CUNY Geometric Analysis
Conference, in 2006 \cite{sz3}.

\section{Technicalities.\label{tech}}
The constructions are performed on 
2-step nilpotent metric Lie groups and their solvable extensions. The
nilpotent groups are defined on 
corresponding orthogonal
direct sums, $\mathbf v\oplus\mathbf z$, of Euclidean spaces where the
components, 
$\mathbf v=\mathbb R^k$ and
$\mathbf z=\mathbb R^l$, 
are called X- and Z-space respectively. The Lie algebra is completely
determined by the linear space, $J_{\mathbf z}$,
of skew endomorphisms whose actions 
on the X-space are defined by the relation 
$\langle [X,Y],Z\rangle =\langle J_Z(X),Y\rangle$, where 
$X,Y\in\mathbf v$
and $J_Z$ is the endomorphism associated with  
$Z\in\mathbf z$. 
The Riemannian metric, $g$, is the left invariant extension of the
natural Euclidean metric on the Lie algebra. The exponential map
identifies the Lie algebra with the group itself, 
thus also the group can be
considered to be defined on the same  
$(X,Z)$-space. Each group, $(N,g)$, 
extends into a solvable group $(SN,g_s)$,
where a point is represented by $(t,X,Z)$. 

Particular 2-step nilpotent Lie groups are the so called
Heisenberg-type groups,
defined by endomorphisms
satisfying the Clifford condition $J^2_Z=-|Z|^2id$. These metric groups 
are attached to Clifford modules, thus the classification
of these modules provides classification also for the H-type groups.
In this case the X-space decomposes into the product 
$\mathbf v=(\mathbf R^{r(l)})^{a+b}=\mathbf R^{r(l)a}\times
\mathbf R^{r(l)b}$ and
endomorphisms $J_Z$ 
are defined by endomorphisms $j_Z$
acting on the smaller space $\mathbf R^{r(l)}$ 
such that they act on  
$\mathbf R^{r(l)a}$ resp.
$\mathbf R^{r(l)b}$ according to the Cartesian product 
$j_Z\times\dots\times j_Z$ resp. $-j_Z\times\dots\times -j_Z$. 
The H-type groups are denoted by the symbol 
$H^{(a,b)}_l$, which indicates the above decomposition. 
The solvable extensions of H-type groups are denoted by $SH^{(a,b)}_l$.  

The Laplacian on a H-type group is of the form
\begin{equation}
\label{Delta}
\Delta=\Delta_X+(1+\frac 1{4}|X|^2)\Delta_Z
+\sum_{\alpha =1}^r\partial_\alpha D_\alpha \bullet,
\end{equation}
where $D_\alpha\bullet$ denotes directional derivative along
the vector field
$X\to J_\alpha (X)=J_{Z_\alpha}(X)$ 
defined for each element,  
$Z_\alpha$, of an orthonormal basis of the Z-space. The integral
curves of this field are the Hopf circles defined for the complex structure 
$J_\alpha$. 
In the isospectrality constructions performed in this paper
one should deal with this compound operator. Earlier, the constructions
were performed on 
center periodic H-type groups,  
$\Gamma\backslash H$, 
defined by factorizing the center of the group
with a Z-lattice $\Gamma =\{Z_\gamma\}$. 
In this case the $L^2$ function
space is the direct sum of function
spaces $W_\gamma$ spanned by functions of the form
$\Psi_\gamma (X,Z)
=\psi (X)e^{2\pi\mathbf i\langle Z_\gamma ,Z\rangle}.
$
Each $W_\gamma$ is invariant under the action of the Laplacian, i. e., 
$\Delta \Psi_{\gamma }(X,Z)=\Box_{\gamma}\psi (X)
e^{2\pi\mathbf i\langle Z_\gamma ,Z\rangle}$,  
where operator $\Box_{\gamma }$, acting on 
$L^2(\mathbf v)$, is of the form
\begin{equation}
\label{Box}
\Box_{\gamma }
=\Delta_X + 2\pi\mathbf i D_{\gamma }\bullet -4\pi^2
|Z_\gamma |^2(1 + 
\frac 1 {4} |X|^2).
\end{equation}
Note that the first operator involves
all endomorphisms $J_\alpha$ while the second one involves,
regarding each invariant subspace $W_\gamma$, 
only $J_\gamma$. 

\noindent{\bf Remark.}
There is pointed out in \cite{sz3,sz4,sz5} that operator (\ref{Box})
is nothing but the Zeeman-Hamilton operator of a free charged
particle (the 2D version is called Landau Hamiltonian), which was
used for explaining the Zeeman effect. Term involving $D_\gamma\bullet$
is the so called angular momentum operator, which represents a
preliminary version of the spin concept. The non-periodic metric
group $(N,g)$ strongly relates to Dirac's relativistic 
multi-time model, which,
in order to furnish relativistic features on the quantum level, endowed
the particles with individual self-times. In the H-type model the
multi-time is represented by the multi-dimensional center of the group.
Regarding this relativistic interpretation, which is not the same as the
classical relativism, Laplacian (\ref{Delta}) on the total space 
(space-time) corresponds to the Klein-Gordon operator. Note, however,
that this multi-time operator is an elliptic one. This fact points to the
distinctive features of the multi-time and classical relativism. Operator,
$\mathbf M$, involving all
angular momentum operators is called compound angular momentum
operator.

\section{Isospectrality constructions.}

The isospectrality constructions are performed on H-type
groups, 
$H^{(a,b)}_l$, 
and on their solvable extensions, 
$SH^{(a,b)}_l$, first. It is only the last section where they are 
extended to 
$\sigma$-equivalent groups whose endomorphism spaces are constructed
by perturbing the Clifford endomorphism spaces. 
The main goal is to describe these constructions on non-periodic groups,
however, in order to see both the similarities and differences, they are 
briefly reviewed here also in the center periodic cases. For fixed $a+b$ 
and $l$, these groups are defined on the same $(X,Z)$- resp. 
$(t,X,Z)$-space. There is established in many different ways that 
the metrics, $g_l^{(a,b)}$, in a family have
completely different local geometries (cf., for instance, the striking 
examples), yet they are 
isospectral on a wide range of submanifolds. 

\subsection{Constructing the ball$\times$torus- and 
sphere$\times$torus-examples.} 
These examples are constructed for a
family,
$\Gamma\backslash H^{(a,b)}_l$, 
of Z-periodic manifolds. The 
submanifolds considered are torus bundles over a ball (resp. sphere) 
around the origin of the X-space.
An intertwining operator can be constructed such that,
for each invariant space $W_\gamma$ constructed above, 
just an orthogonal transformation conjugating 
$J_{\gamma}^{(a,b)}$ to
$J_{\gamma}^{(a^\prime,b^\prime)}$ on the X-space should be considered.
The intertwining operator on $W_\gamma$ is defined by the map induced
on functions $\psi (X)$, defined in the Fourier-Weierstrass decomposition,
by this point transformation.
This transformation clearly
intertwines $\Box_\gamma$ with $\Box^\prime_\gamma$ such that it
keeps also the boundary conditions. 
(The boundary conditions can be described
in terms of radial functions. The intertwining 
of boundary conditions is then  
due to the invariance of these functions under the action
of the operator.) It induces an appropriate
intertwining operator also on the boundary manifolds.
The striking examples appear on the quaternionic families 
$H^{(a,b)}_3$, in which case the 
sphere$\times$torus-type boundary manifolds in 
$\Gamma\backslash H^{(a+b,0)}_3$
are homogeneous while the others in the family are
locally inhomogeneous. Note that the simplicity of this case
is due to the fact that the intertwining operator is constructed, on each 
invariant space $W_\gamma$ separately, 
by a single endomorphism, $J_\gamma$. 
 
\subsection{The ball- and sphere-type domains.} These
examples were originally constructed in \cite{sz1,sz2}.
The ball-type domains are, by definition,  
diffeomorphic to Euclidean balls
such that the sphere-type boundary manifolds
are level sets described by equations of the form 
$\varphi (|X|,|Z|)=0$ resp. $\varphi (|X|,|Z|,t)=0$. 
These domains are invariant under the action of the orthogonal
group $\mathbf O (\R^k)\times\mathbf O (\R^l)$, thus, they may be
called domains of 
$(X,Z)$-revolutions. They can be visualized such that
there is an X-ball centered at the origin of the X-space considered 
over the 
points of which
there are Z-balls of radius $R_Z(|X|)$ around the origin
of the Z-space considered. Then, the boundary is a level set defined by
the equation
$\varphi (|X|,|Z|)=|Z|-R_Z(|X|)=0$. By this reason, function $\varphi$ is
called radius function.

Note that radius $R_Z(|X|)$
is constant along a sphere defined by a constant radius $R_X=|X|$
in the X-space. Furthermore, the
ball bundle defined by the Z-balls 
over this X-sphere is trivial. These are
the so called {\it sphere$\times$ball-type manifolds} whose boundaries
are {\it sphere$\times$sphere-type manifolds}. The new examples,
not discussed in the earlier papers, are constructed on these
domains and surfaces.

An other visualization
can be started out with a Z-ball in the Z-space over the points of
which there are X-balls of radius $R_X(|Z|)$ considered.
Then the boundary is defined by $|X|-R_X(|Z|)=0$. However, this paper
proceeds with the first description.

In the solvable
case one should consider $(Z,t)$-balls and $(Z,t)$-spheres around 
the origin $(0_Z,1)$ defined for the hyperbolic $(Z,t)$-space. 
The base manifold is the same X-sphere as before.
Note that a Z-ball $B_{R_Z}(0_Z)$ (resp. Z-sphere $S_{R_Z}(0_Z)$)
uniquely extends into a geodesic ball (resp. sphere) of the  
hyperbolic $(Z,t)$-space. 
A sphere$\times$ball-type domain can be described as a hypersurface
in a ball-type domain such that the Z-balls (resp. $(Z,t)$-balls) of the
ball type domain are considered only over the points of a sphere $S_{R_X}$ 
laying in the X-space. 
Similarly, the sphere$\times$sphere-type manifolds can be regarded
as hypersurfaces in the sphere-type manifolds.

The isospectrality will be
investigated, first, for the discrete families, 
$H^{(a,b)}_l$, defined by the same $a+b$ and $l$. 
The Laplacian  
is described then by (\ref{Delta}). Comparing with the Zeeman operator
(\ref{Box}), this operator involves all the endomorphisms, making
the constructions much more difficult. The Laplacians of the members in
a family differ from each other just by the last term, $\mathbf M$,
which is called compound  angular momentum operator. The spectral
investigation both of $\mathbf M$ and $\Delta$ is completely missing
in the literature. 
Note that this most intriguing operator, $\mathbf M$,  
commutes with both operators in the rest part of (\ref{Delta}). 

\subsection{Eigenfunctions motivating the intertwining operators.} 
The eigenfunctions constructed next
are not directly used in the isospectrality constructions and 
the rest part of this paper is understandable without knowing about 
their actual
explicit form described in the second half of this section. However,
there are important concepts introduced in the first part which are
heavily used later on.
The ultimate reason for describing these functions here is that they
very clearly suggest the explicit 
form of the sought intertwining operators.
 
Since $\mathbf M$ commutes with the rest part, $\mathbf O$, 
of $\Delta$, the eigenfunctions of $\Delta$ can be
sought as common eigenfunctions for 
both operators $\mathbf M$ and $\mathbf O$.
In the very first step we look for the eigenfunctions of a single 
angular momentum operator $\mathbf D_V\bullet$, defined by a
Z-vector $V$.  For a fixed X-vector $Q$
and unit Z-vector $V_u={1\over |V|}V$, consider the X-function
$\Theta_{Q}(X,V_u)=\langle Q+\mathbf iJ_{V_u}(Q),X\rangle$ and its 
conjugate $\overline{\Theta}_{Q}(X,V_u)$. For vector $V=|V|V_u$,
these functions are 
eigenfunctions of $D_{V}\bullet$ with eigenvalue $-|V|\mathbf i$ resp.
$|V|\mathbf i$. The higher order eigenfunctions are of the form
$\Theta_{Q}^p\overline\Theta^q_{Q}$ 
with eigenvalue $(q-p)|V|\mathbf i$.
 
In order to find eigenfunctions of the compound operator $\mathbf M$,
consider a sphere $S_{R_Z}$ of radius $R_Z$
around the origin in the Z-space. For an appropriate function
$\phi (|X|,V)$, depending on $|X|$ and $V\in S_{R_Z}$,
define 
\begin{equation}
\label{FourR_Z}
\mathcal F_{QpqR_Z}(\phi )(X,Z)
=\oint_{S_{R_Z}}e^{\mathbf i
\langle Z,V\rangle}\phi (|X|,V)
\Theta_{Q}^p(X,V_u)\overline\Theta^q_{Q}(X,V_u)dV.
\end{equation} 
By $V_u=V/|V|$, the $V_u$ is considered as a function depending on $V$. 
Due to the relation $\mathbf M\oint =\oint \mathbf iD_V\bullet$, 
this function is an eigenfunction
of $\mathbf M$ with the real eigenvalue $(p-q)R_Z$. 
These functions are eigenfunctions also of $\Delta_Z$ with eigenvalue
$R_Z^2$. Also note 
that these eigenvalues do not change by varying $Q$, or, if the
simple functions
$
\Theta_{Q}^p(X,V_u)
$
resp.
$
\overline\Theta^q_{Q}(X,V_u)
$
are exchanged for their pluralistic versions
$
\Theta_{Q_{1\dots p}}(X,V_u):=
\Theta_{Q_1}(X,V_u)\dots
\Theta_{Q_p}(X,V_u)
$
resp.
$
\overline\Theta_{Q_{p+1\dots p+q}}(X,V_u)
$.

Functions (\ref{FourR_Z}) 
defined by simple resp. pluralistic 
functions are said to be one-pole resp.
multiple-pole functions with poles $Q$ resp. $\{Q_i\}$. 
The function space generated for fixed 
1-pole (resp. multi-poles) by all 
possible $\phi$ is not invariant with respect to 
the action of $\Delta_X$, thus the eigenfunctions of the complete
operator $\Delta$ do not appear in this form. In order to find the
common eigenfunctions,
the homogeneous but non-harmonic 1-pole polynomials 
$
\Theta_{Q}^p\overline\Theta^q_{Q}
$
(resp. the multiple-pole polynomials)
of the X-variable should be exchanged for the 
1-pole harmonic polynomials
$
\Pi^{(n)}_X (\Theta_{Q}^p\overline\Theta^q_{Q})
$, 
(resp. to the corresponding multiple-pole harmonic polynomials)
defined by the projection, $\Pi^{(n)}_X$, 
onto the space of $n=(p+q)$-order 
homogeneous harmonic 
polynomials of the X-variable. These projections 
are explicitly described in the form 
\begin{equation}
\label{pi_x}
\Pi^{(n)}_X =\Delta_X^0+
B^{(n)}_1
|X|^2\Delta_X +B^{(n)}_2|X|^4\Delta_X^2+\dots
\end{equation}
in \cite{sz2} (cf. formula (3.14) there), 
where $\Delta_X^0=id$ and the constants  
$B^{(n)}_i$
are determined by a recursion formula. 
This formula easily implies that also  
\begin{equation}
\label{HFourR_Z}
\mathcal {HF}_{QpqR_Z}(\phi )(X,Z)
=\oint_{S_{R_Z}}e^{\mathbf i
\langle Z,V\rangle}\phi (|X|,V)
\Pi_X(\Theta_{Q}^p(X,V_u)\overline\Theta^q_{Q}(X,V_u))dV
\end{equation}
are eigenfunctions
of $\mathbf M$ and $\Delta_Z$ with eigenvalues 
belonging to (\ref{FourR_Z}). The same statement is true regarding the 
multiple-pole-cases. 

The action of the complete Laplacian (\ref{Delta}) is a combination of
X-radial differentiation, $\partial_{|X|}$, 
and multiplications with functions depending 
just on $|X|$. I. e., the action is completely reduced to X-radial
functions. Also this reduced form of the Laplacian is not changing by 
varying $Q$, or, switching to multiple-pole functions. 
The eigenfunctions of $\Delta$ can be found by seeking
the eigenfunctions of the reduced operator among the X-radial functions.
The explicit computations are carried out in \cite{sz4,sz5}. Since these 
details are not used in this paper, we just indicate that the
eigenfunctions appear in the form
$
\oint_{S_{R_Z}}e^{\mathbf i
\langle Z,V\rangle}\phi (V)
F^{(p,q)}_{Q}(X,V_u))dV,
$
(resp. in a corresponding multiple-pole version of this function),
where
$
F^{(p,q)}_{Q}(X,V_u))
$
is an eigenfunction of operator $\Box_\gamma$ satisfying 
$|V_\gamma |=|V|=R_Z$. The latter ones are explicitly described in 
\cite{sz4,sz5} in terms of homogeneous harmonic polynomials which are
multiplied with radial functions. The eigenfunctions are determined
below also by a different method, using the so called It\^o polynomials.

Note that this construction is carried out for a fixed 1-pole $Q$ (resp. a 
fixed multiple-pole,
$\{Q_1,\dots ,Q_p,Q_{p+1},\dots ,Q_{p+q}\}$).
An other type of constructions is as follows.
For any unit vector $V_u$ of the Z-space, 
consider a complex orthonormal basis 
$\{Q_{V_u1},\dots ,Q_{V_uk/2}\}$
on the complex X-space defined by the complex structure $J_{V_u}$
such that the vectors in front lay in 
$\mathbf v^{(a)}$ 
and all the others are in 
$\mathbf v^{(b)}$.
Such a basis defines the complex coordinate system
$\{z_{V_u1}=\Theta_{Q_{1V_u}},\dots ,z_{V_uk/2}=\Theta_{Q_{(k/2)V_u}}\}$ 
on the X-space. This basis field must be smooth on an everywhere 
dense open subset of
the unit Z-sphere such that it 
is the complement of a set of $0$ measure. For given values
$p_1,q_1,\dots ,p_{k/2},q_{k/2}$,  
consider the polynomial
$
\prod_{i=1}^{k/2}
z_{V_ui}^{p_i} 
\overline z_{V_ui}^{q_i}. 
$
Then the functions  
\begin{equation}
\label{Fourprod}
\oint_{S_{R_Z}}e^{\mathbf i
\langle Z,V\rangle}\phi (V)\varphi (|X|)
\prod_{i=1}^{k/2}
z_{V_ui}^{p_i} 
\overline z_{V_ui}^{q_i} dV
\end{equation}
are eigenfunctions of the compound angular momentum operator $\mathbf M$.
In order to have an eigenfunction for the complete Laplacian, one can
use the above described method of projecting the polynomial into the 
space of homogeneous harmonic polynomials which have order
$(p_1+q_1+\dots +p_{k/2}+q_{k/2})$. 
In \cite{sz4,sz5}, the 
eigenfunctions of $\Box_\gamma$ are determined also by an other 
method, seeking them in the form
$
 h^{(p_iq_i)}_{V_u}e^{-R_Z\langle z_{V_u},\overline z_{V_u}\rangle}, 
$
where the
$
 h^{(p_iq_i)}_{V_u} 
$ is an
$(p_1+q_1+\dots +p_{k/2}+q_{k/2})$-order polynomial. Then this
function is an eigenfunction of $\Box_\gamma$ satisfying 
$|Z_\gamma |=R_Z$ if and only if the latter function is an It\^o
polynomial regarding the complex structure $J_{V_u}$. The final form of
the eigenfunction is
\begin{equation}
\oint_{S_{R_Z}} e^{\mathbf i\langle Z,V\rangle}
\phi (V)
 h^{(p_iq_i)}_{V_u} 
e^{-R_Z z_{V_u}\cdot\overline z_{V_u}}
dV.
\end{equation}
Since It\^o's polynomials are non-homogeneous, this is a different
representation of the eigenfunctions. These explicit descriptions
of the eigenfunctions will not be used in the rest part of the paper.

\section{Constructing the intertwining operators.\label{io}}
\subsection{Constructions on ball-type domains.}

The constructions described in this
sections are carried out
for Heisenberg type groups,
$H^{(a,b)}_l$ and $H^{(a^\prime,b^\prime)}_l$, which are in the same 
isospectrality family, i.e., 
$a+b=a^\prime+b^\prime$ holds.  
From each eigenfunction-construction
described above one can derive the corresponding
intertwining operator intertwining the corresponding eigenfunctions
provided by the construction. Note that the functions appearing there are
not of class $L^2$ regarding the Z-variable thus integral formulas
(\ref{FourR_Z}), 
(\ref{HFourR_Z}), 
(\ref{Fourprod}) 
can not be directly used for defining the operator. This is why function 
$\phi (|X|,V)$, depending on $|X|$ and $V\in S_{R_Z}$ originally, is
exchanged for one which an $L^2$ function of the V-variable, 
for any fixed $|X|$,
and the integral is taken over the whole Z-space $\mathbb R^l$. In other
words, the Z-Fourier transform on $L_Z^2$-setting is considered.

Also the order for introducing the various versions of the intertwining
operators is 
an important issue.
The first version is defined for a fixed basis, 
$
\mathbf Q_F
=\{Q_{1},\dots ,Q_{k/2}\}
=\{\mathbf Q^{(a)}_F,
\mathbf Q^{(b)}_F\}
$, 
which does not depend on $V_u$, where the first 
$k^{(a)}/2$ number of vectors 
are in $\mathbf v^{(a)}$ and 
the following $k^{(b)}/2$ number of vectors
are in $\mathbf v^{(b)}$. If these vectors are
chosen such that they form an
orthonormal basis regarding a fixed $J_{V_{0u}}$, then they
form a complex (in general, non-orthonormal) basis for $V_u$'s which
form an everywhere dense open subset on the unit sphere of the
Z-space. This operator is defined by means of the polynomials written
up in terms of the coordinate functions 
\begin{equation}
\label{z(X)}
z_{V_u1}(X)
=\Theta_{Q_1}(X,V_u),
\dots ,z_{V_uk/2}(X)=\Theta_{Q_{k/2}}(X,V_u)
\end{equation}
resp. 
$\{z_{V_u1}^\prime ,\dots ,z_{V_uk/2}^\prime\}$. The denotation indicates
that, although the basis is
fixed, these coordinate functions depend on $V_u$. 

For constructing the eigenfunctions in the previous section, 
one is using a fixed Z-sphere and 
functions, $\phi (V)$, of Dirac-type concentrated on this sphere. Then
the eigenfunctions are represented by Z-Fourier transforms of these
Dirac-type functions. Whereas, in defining the intertwining operator 
$\kappa_{\mathbf Q_F}$ determined for a constant basis $\mathbf Q_F$, one 
considers Z-Fourier
transforms of appropriate $L_Z^2$-functions. More precisely, for an 
$L_Z^2$-function, $\phi (V)$, defined on the Z-space and 
arbitrary set
$\{p_i,q_i\}$, where $i=1,\dots ,k/2$, of natural numbers 
the intertwining operator is defined by 
\begin{eqnarray}
\label{io1}
\kappa_{\mathbf Q_F} : &
\mathcal F_{\mathbf Q_F\{p_iq_i\}}(\phi )(X,Z)=
\int_{\mathbb R^l} e^{\mathbf i\langle Z,V\rangle}
\phi (V)
\prod_{i=1}^{k/2}z^{p_i}_{V_ui}(X)\overline z^{q_i}_{V_ui}(X)dV 
\\
& \to \mathcal F^\prime_{\mathbf Q_F\{p_iq_i\}}(\phi )(X,Z)=
\int_{\mathbb R^l} e^{\mathbf i\langle Z,V\rangle}
\phi (V)
\prod_{i=1}^{k/2}z^{\prime p_i}_{V_ui}(X)
\overline z^{\prime q_i}_{V_ui}(X)dV.
\nonumber 
\end{eqnarray}
I. e., the $\kappa_{\mathbf Q_F}$ corresponds to a function, 
which is defined by
the Z-Fourier transform formula in terms of
$\phi , z_i^{p_i},\overline z_i^{q_i}$,
the function defined by the same expression
but which is written up in terms of 
$\phi , z_i^{\prime p_i}$, and $\overline z_i^{\prime q_i}$.
In these formulas, the $V_u=V/|V|$ is considered as a function 
depending on $V$, furthermore, the dependence of the complex coordinate
functions on the X-variable is described 
in (\ref{z(X)}). Note that in this
first version of the operator function $\phi$ depends just on $V$ and 
not on $|X|$.

The domain and range of this operator is discussed in the next section.
In this section one is focusing on the well-definedness and the 
intertwining of the Laplacians. Concerning these questions, we have.

\begin{theorem}
The above 
$\kappa_{\mathbf Q_F}$ 
is a well defined one-to-one operator.
\end{theorem}
\begin{proof}
This theorem is well established
by proving that the image of a function which is in the domain of 
$\kappa_{\mathbf Q_F}$
and vanishes almost everywhere is a function vanishing almost everywhere. 
For proving this statement,
suppose that function
\[
\tilde\varphi (X,Z)=
\int_{\mathbb R^l} e^{\mathbf i\langle Z,V\rangle}
\sum_
{\{p_i,q_i\}}
\phi_ 
{\{p_i,q_i\}}
(V)
\prod_{i=1}^{k/2}z^{p_i}_{V_ui}(X)\overline z^{q_i}_{V_ui}(X),
\] 
where the terms of the sum (series) are defined regarding the independent
polynomials
$
\prod_{i=1}^{k/2}z^{p_i}_{V_ui}\overline z^{q_i}_{V_ui},
$
vanishes almost everywhere. Since the Z-Fourier transform is a one-to-one
map on the corresponding $L^2$-Hilbert space, for any fixed $X$, function
$
\varphi (X,V) =\sum_
{\{p_i,q_i\}}
\phi_ {\{p_i,q_i\}}
(V)
\prod_{i=1}^{k/2}z^{p_i}_{V_ui}(X)\overline z^{q_i}_{V_ui}(X)
$ 
(whose Fourier transform is considered) must vanish almost everywhere.
By the independence of the polynomials
$
\prod_{i=1}^{k/2}z^{p_i}_{V_ui}\overline z^{q_i}_{V_ui},
$ 
what is satisfied for almost all $V_u$, 
a general function $\varphi$ defined by this formula is non-zero 
if there is
a non-zero $L^2$-function
$
\phi_ {\{ p_i,q_i\}}
$ among the component functions.
Therefore, due to the assumption of this theorem, all these $\phi$'s 
must vanish almost everywhere. Then, also the image 
$\tilde\varphi^\prime$ must
vanish almost everywhere. This proves
the statement completely.

Actually, one has proved the following stronger statement: Function
$\tilde\varphi^\prime$ 
is zero almost everywhere if and only if its preimage
$\tilde\varphi$ 
is zero almost everywhere. Thus also the 
one-to-one property is established
completely.
\end{proof}

Much more handy alternative definitions of the very same 
$\kappa_{\mathbf Q_F}$ are established in the following theorem.

\begin{theorem}
(A) In the definition of 
$\kappa_{\mathbf Q_F}$, 
function $\phi$ may depend
also on $|X|$, or, even more,on 
$
|X^{(a)}|
$
and 
$
|X^{(b)}|
$, where
$X=X^{(a)}+X^{(b)}$ corresponds to the decomposition 
$\mathbf v=\mathbf v^{(a)}\oplus\mathbf v^{(b)}$. 
I. e., an equivalent version is: 
\begin{eqnarray}
\label{io2}
\kappa_{\mathbf Q_F} : 
\mathcal F_{\mathbf Q_F\{p_iq_i\}}(\phi )(X,Z)=
\int_{\mathbb R^l} e^{\mathbf i\langle Z,V\rangle}
\phi (|X|,V)
\prod_{i=1}^{k/2}z^{p_i}_{V_ui}(X)\overline z^{q_i}_{V_ui}(X)dV \\
\to \mathcal F^\prime_{\mathbf Q_F\{p_iq_i\}}(\phi )(X,Z)=
\int_{\mathbb R^l} e^{\mathbf i\langle Z,V\rangle}
\phi (|X|,V)
\prod_{i=1}^{k/2}z^{\prime p_i}_{V_ui}(X)
\overline z^{\prime q_i}_{V_ui}(X)dV.
\nonumber 
\end{eqnarray}
In the more general version, the 
$\phi (|X|,V)$
is replaced by
$
\phi (|X^{(a)}|,|X^{(b)}|,V).
$

(B) The 
$\kappa_{\mathbf Q_F}$ 
intertwines both the Euclidean Laplacian $\Delta_X$ and the projections
$\Pi^{(n)}_X$, 
$\Pi^{(n_a)}_X$, 
$\Pi^{(n_b)}_X$ with themselfs respectively. These projections are  
described in (\ref{pi_x}), furthermore, 
$n=p+q=\sum p_i+\sum q_i$,
$n_a=p_a+q_a=\sum_{i=1}^a (p_i+ q_i)$,
$n_b=p_b+q_b=\sum_{i=a+1}^{a+b} (p_i+ q_i)$.
The operator can be written 
in the following alternative form:
\begin{eqnarray}
\label{io3}
\kappa_{\mathbf Q_F} : 
\mathcal{HF}_{\mathbf Q_F\{p_iq_i\}}(\phi )(X,Z)=
\int_{\mathbb R^l} e^{\mathbf i\langle Z,V\rangle}
\phi (|X|,V)
\Pi^{(n)}_X(
\prod_{i=1}^{k/2}z^{p_i}_{V_ui}\overline z^{q_i}_{V_ui})dV \\
 \to \mathcal{HF}^\prime_{\mathbf Q_F\{p_iq_i\}}(\phi )(X,Z)=
\int_{\mathbb R^l} e^{\mathbf i\langle Z,V\rangle}
\phi (|X|,V)\Pi^{(n)}_X(
\prod_{i=1}^{k/2}z^{\prime p_i}_{V_ui}\overline z^{\prime q_i}_{V_ui})dV.
\nonumber 
\end{eqnarray}
In the more general version, the 
$\phi (|X|,V)$ resp.
$
\Pi^{(n)}_X(\prod_{i=1}^{k/2}\dots )
$
are replaced by
$
\phi (|X^{(a)}|,|X^{(b)}|,V) 
$
resp.
$
\Pi^{(n_a)}_X(\prod_{i=1}^{a}\dots ) 
\Pi^{(n_b)}_X(\prod_{i=a+1}^{a+b}\dots ).
$

(C) These versions for defining
$\kappa_{\mathbf Q_F}$
allow to introduce its domain in a more precise way.  
In order to work on $L^2$ function spaces, one should
consider functions of the form 
$\phi (|X|, V)=e^{-|X|^2}\varphi (|X|,V)$, 
where $\varphi$, 
for any fixed $|X|$, 
is of class $L^2$ with respect to the V-variable,
and it
is a polynomial with respect to the $|X|$-variable. By plugging them into
the Fourier transform formula, they 
generate a pre-Hilbert space whose closure, regarding the $L^2$-Hilbert
norm, is a Hilbert space. (In the
next section, this domain is identified with the standard $L^2$-Hilbert 
space defined on $\mathbb R^k\oplus\mathbb R^l$.) A larger domain can be
generated by functions 
$\phi (|X|,V)$, 
where, keeping the above assumption
regarding the $V$-variable, function  
$\phi (|X|,\dot{V})$ depending on variable $|X|$ is of class $L^2$ 
for almost all fixed $\dot{V}$ on
any interval $0\leq |X|\leq R$. 
\end{theorem}
\begin{proof}
(A) This proof explores that the system 
$\mathbf Q_F$ of independent vectors decomposes into two subsystems,
$\mathbf Q_F=\{\mathbf Q^{(a)},\mathbf Q^{(b)}\}$,
consisting vectors from 
$\mathbf v^{(a)}$
resp.
$\mathbf v^{(b)}$. They form a complex basis, both for 
$J_{V_u}$
and
$J^\prime_{V_u}$,
for an everywhere dense open set of unit vectors $V_u$. This basis can
not be orthonormal for all $V_u$, even though it is orthonormal for some
$V_u$.
For each
$V_u$, let
$\mathbf R_{V_u}=\{\mathbf R_{V_u}^{(a)},\mathbf R_{V_u}^{(b)}\}$
be a complex orthonormal basis, regarding the complex structure
$J_{V_u}$,
defining the complex coordinate system
$\{z_{\mathbf R_{V_u}i}\}$.
Then, there exist complex matrix, 
$
(c_{ij}(V_u))
$, such that 
$
z_{\mathbf R_{V_u}i}
=\sum_jc_{ij}(V_u)
z_{\mathbf QV_uj}
$ 
hold. The equation expressing the orthonormality is 
$
\sum_{mn}
c_{im}
\overline c_{jn}
\langle Q_m,Q_n\rangle
=\delta_{ij},
$ 
where matrix with entries
$
\langle Q_m,Q_n\rangle
$
is real. The basis field can be chosen such that it is continuous on an 
everywhere dense open subset of the unit vectors $V_u$.

In terms of $J^\prime_{V_u}$,
the above matrix transformation 
defines an other complex coordinate system,
$
z^\prime_{
\mathbf R^\prime_{V_u}i}
=\sum_jc_{ij}(V_u)
z^\prime_{\mathbf Q_FV_uj}
$,
with the corresponding complex basis
$
\mathbf R^\prime_{V_u}
$. 
Although it is not the same as 
$
\mathbf R_{V_u}
$, but, 
due to the relations
$J^{\prime (a)}_{V_u}=J^{(a)}_{V_u}$ and
$J^{\prime (b)}_{V_u}=-J^{(b)}_{V_u}$,
this new basis also yields the above orthonormality equation. 
Therefore, it 
is orthonormal regarding the complex structure
$J^\prime_{V_u}$ and, according to the computations: 
\begin{eqnarray}
|X|^2=
\sum_i 
z_{\mathbf R_{V_u}i}
\overline z_{\mathbf R_{V_u}i}=
\sum_i
(\sum_jc_{ij}(V_u)
z_{\mathbf Q_FV_uj})
(\sum_r\overline c_{ir}(V_u)
\overline z_{\mathbf Q_FV_ur})=
\end{eqnarray}
\begin{eqnarray}
=\sum_{j,r}
(\sum_ic_{ij}(V_u)
\overline c_{ir}(V_u))
z_{\mathbf Q_FV_uj}
\overline z_{\mathbf Q_FV_ur})
\nonumber
\end{eqnarray}
\begin{eqnarray}
=\sum_i 
z^\prime_{\mathbf R_{V_u}i}
\overline z^\prime_{\mathbf R_{V_u}i}=
\sum_{j,r}
(\sum_ic_{ij}(V_u)
\overline c_{ir}(V_u))
z^\prime_{\mathbf Q_FV_uj}
\overline z^\prime_{\mathbf Q_FV_ur})
\nonumber
\end{eqnarray}
the X-radial function $|X|^2$ appears in 
the same polynomial form
regarding both coordinate systems
$
\{z_{\mathbf Q_FV_uj}\}
$
and
$
\{z^\prime_{\mathbf Q_FV_uj}\}
$. The same statement is true for any power, or, by using power series, 
for any function 
$\phi (|X|)$ of the basic radial function. This observation proves (A)
completely.

(B) Regarding the Laplacian the same computation yield:
\begin{eqnarray}
\Delta_X=
\sum_i 
\partial_{z_{\mathbf R_{V_u}i}}
\partial_{\overline z_{\mathbf R_{V_u}i}}=
\sum_i
(\sum_jc_{ij}(V_u)
\partial_{z_{\mathbf Q_FV_uj}})
(\sum_r\overline c_{ir}(V_u)
\partial_{\overline z_{\mathbf Q_FV_ur}})=\\
=\sum_{j,r}
(\sum_ic_{ij}(V_u)
\overline c_{ir}(V_u))
\partial_{z_{\mathbf Q_FV_uj}}
\partial_{\overline z_{\mathbf Q_FV_ur})}
\nonumber\\
=\sum_i 
\partial_{z^\prime_{\mathbf R_{V_u}i}}
\partial_{\overline z^\prime_{\mathbf R_{V_u}i}}=
\sum_{j,r}
(\sum_ic_{ij}(V_u)
\overline c_{ir}(V_u))
\partial_{z^\prime_{\mathbf Q_FV_uj}}
\partial_{\overline z^\prime_{\mathbf Q_FV_ur}}),
\nonumber
\end{eqnarray}
i. e., also this Laplacian appears regarding both coordinate system in
the same form. This proves the invariance 
of $\Delta_X$ under the action of 
$\kappa_{\mathbf Q_F}$. The explicit formula (\ref{pi_x}) regarding 
$\Pi^{(n)}_X$ along with (A) and
the above statement concerning $\Delta_X$ prove (B) completely.

(C) This statement is self-contained.
\end{proof}

Yet an other alternative definition of
$\kappa_{\mathbf Q_F}$
can be introduced by using 1-pole functions
$\Theta_{Q}^p\overline\Theta^q_{Q}$,
where $Q$ is in the real span of the vector-system $\mathbf Q_F$, i. e., 
$Q\in Span_{\mathbb R}(\mathbf Q_F)$.
The version using multi-pole functions, defined by vectors laying in
$Span_{\mathbb R}(\mathbf Q_F)$,
adds nothing new to the 
above polynomial-version, therefore, this case is omitted here. Also 
note that the multi-pole functions span the same function space spanned 
by the 1-pole functions, thus, complete analysis can be performed
by using only the simple ones. 
\begin{theorem} The 
$\kappa_{\mathbf Q_F}$ is well defined by each of the following 
correspondences:
\begin{eqnarray}
\label{io4}
\kappa_{\mathbf Q_F } : &
\mathcal F_{Qpq}(\phi )(X,Z)=
\int_{\mathbb R^l} e^{\mathbf i\langle Z,V\rangle}
\phi (V)
\Theta^{p}_{QV_u}(X)\overline \Theta^{q}_{QV_u}(X)dV 
\\
& \to \mathcal F^\prime_{Qpq}(\phi )(X,Z)=
\int_{\mathbb R^l} e^{\mathbf i\langle Z,V\rangle}
\phi (V)
\Theta^{\prime p}_{QV_u}(X)
\overline \Theta^{\prime q}_{QV_u}(X)dV,
\nonumber 
\end{eqnarray}
\begin{eqnarray}
\label{io5}
\kappa_{\mathbf Q_F } : &
\mathcal F_{Qpq}(\phi )(X,Z)=
\int_{\mathbb R^l} e^{\mathbf i\langle Z,V\rangle}
\phi (|X|,V)
\Theta^{p}_{QV_u}(X)\overline \Theta^{q}_{QV_u}(X)dV 
\\
& \to \mathcal F^\prime_{Qpq}(\phi )(X,Z)=
\int_{\mathbb R^l} e^{\mathbf i\langle Z,V\rangle}
\phi (|X|,V)
\Theta^{\prime p}_{QV_u}(X)
\overline \Theta^{\prime q}_{QV_u}(X)dV,
\nonumber 
\end{eqnarray}
\begin{eqnarray}
\label{io6}
\kappa_{\mathbf Q_F } : &
\mathcal{HF}_{Qpq}(\phi )(X,Z)=
\int_{\mathbb R^l} e^{\mathbf i\langle Z,V\rangle}
\phi (|X|,V)
\Pi^{(p+q)}_X
(\Theta^{p}_{QV_u}\overline \Theta^{q}_{QV_u})dV 
\\
& \to \mathcal{HF}^\prime_{Qpq}(\phi )(X,Z)=
\int_{\mathbb R^l} e^{\mathbf i\langle Z,V\rangle}
\phi (|X|,V)
\Pi^{(p+q)}_X(
\Theta^{\prime p}_{QV_u}
\overline \Theta^{\prime q}_{QV_u})dV,
\nonumber 
\end{eqnarray}
where 
$Q\in Span_{\mathbb R}(\mathbf Q_F)$ is an arbitrary vector and
$
\Theta_{QV_u}(X)
:=\Theta_{Q}(X,V_u).
$

There is a more general version also in 
this case which corresponds to the
exchange of $Q$ for a pair, 
$(Q^{(a)},Q^{(b)})$, 
followed by the exchange of functions 
$\phi (|X|,V)$,\,
$
\Theta^{p}_{QV_u}\overline \Theta^{q}_{QV_u}
$,\,
$
\Pi^{(p+q)}_X(..)
$
for the following ones
\[
\phi (|X^{(a)}|,|X^{(b)}|,V),\,\,
\Theta^{p_a}_{Q^{(a)}V_u}\overline \Theta^{q_a}_{Q^{(a)}V_u} 
\Theta^{p_b}_{Q^{(b)}V_u}\overline \Theta^{q_b}_{Q^{(b)}V_u},\,\,
\Pi^{(p_a+q_a)}_{X^{(a)}}(..) 
\Pi^{(p_b+q_b)}_{X^{(b)}}(..)
\]
respectively.
\end{theorem}

\noindent{\bf Remark.} Although, it is not 
defined by an appropriate basis,
the well-defined\-ness of the operator is not 
jeopardized in this theorem. It can be defined, however, by constructing
a basis as follows.

For fixed $Q$ and natural numbers $p$ and $q$, let
$
\mathbf \Phi_{Qpq}
$
(resp. 
$
\mathbf \Xi_{Qpq}
$)
be the $L^2$ function space spanned by functions of the form
$\mathcal F_{Qpq}(\phi )(X,Z)$
(resp. $\mathcal{HF}_{Qpq}(\phi )(X,Z)$),
where $\phi(|X|,V)$
can be an arbitrary $L_Z^2$-function. For fixed $p$ and $q$, all these 
spaces sum up to the total spaces
$
\mathbf \Phi_{pq}=\sum_Q
\mathbf \Phi_{Qpq}
$
(resp.
$
\mathbf \Xi_{pq}=\sum_Q
\mathbf \Xi_{Qpq}
$). There exist finite many $Q_i$  such
that the total space is the direct sum of the independent subspaces 
$
\mathbf \Phi_{Q_ipq}
$
(resp. 
$
\mathbf \Xi_{Q_ipq}
$). For the two type of total spaces these numbers are different. Due to
the non-degeneracy of
$\kappa_{\mathbf Q_F}$, 
the total spaces
$
\mathbf \Phi^\prime_{pq}
$
(resp.
$
\mathbf \Xi^\prime_{pq}
$)
are the direct sums of the independent subspaces 
$
\mathbf \Phi^\prime_{Q_ipq}
$
(resp. 
$
\mathbf \Xi^\prime_{Q_ipq}
$). The  
$\kappa_{\mathbf Q_F}$ 
can be defined just by its actions 
$
\kappa_{\mathbf Q_F}:
\mathbf \Phi_{Q_ipq}\to
\mathbf \Phi^\prime_{Q_ipq}
$
(resp. 
$
\kappa_{\mathbf Q_F}:
\mathbf \Xi_{Q_ipq}\to
\mathbf \Xi^\prime_{Q_ipq}
$) on these subspaces. 
This construction method, 
which will not be used in this paper,
can be applied for explicit spectral computations. 
\qed  

Above, six versions of the very same intertwining operator defined by a
constant complex basis, 
$\mathbf Q_F$,
were introduced. 
A changing
orthonormal complex basis field,
$\mathbf Q(V_u)=\{\mathbf Q^{(a)}(V_u),\mathbf Q^{(b)}(V_u)\}$,
which is supposed to be continuous on an everywhere dense open subset 
of the unit vectors
$V_u$,
also defines an intertwining operator. 
Unlike the constant field, which can be orthonormal only for
$V_u$'s of zero measure, the changing field is supposed to be orthonormal
almost everywhere. Then, by 
$J^{\prime (a)}_{V_u}=J^{(a)}_{V_u}$ and
$J^{\prime (b)}_{V_u}=-J^{(b)}_{V_u}$,
also the basis 
$
\mathbf Q^\prime (V_u)
$
is orthonormal regarding the complex structure
$J^\prime_{V_u}$. 
The starting version of the intertwining operator
$\kappa_{\mathbf Q(V_u)}$
defined by a changing complex orthonormal basis is introduced by 
formula (\ref{io1}), where the complex coordinates, $z_{V_ui}$,
are defined by  
$\mathbf Q(V_u)$.
Now we have
\begin{theorem}
\label{euclapl}
Operator
\begin{equation}
\kappa_{\mathbf Q(V_u)}:\quad  
\mathcal F_{Q(V_u)\{p_i,q_i\}}(\phi )(X,Z)\to
\mathcal {F}^{\prime}_{Q(V_u)\{p_i,q_i\}}(\phi )(X,Z),
\end{equation}
defined for all sets, $\{p_i,q_i\}$, of natural numbers and 
$L^2$-functions $\phi (V)$, is a well-defined one-to-one map
leaving the X-radial functions, the Laplacian $\Delta_X$ and the
projections $\Pi_X^{(n)}$ invariant. Thus, this operator can be
defined in the alternative ways 
\begin{equation}
\kappa_{\mathbf Q(V_u)}:\quad  
\mathcal F_{Q(V_u)\{p_i,q_i\}}(\phi )(X,Z)\to
\mathcal {F}^{\prime}_{Q(V_u)\{p_i,q_i\}}(\phi )(X,Z),
\end{equation}
\begin{equation}
\kappa_{\mathbf Q(V_u)}:\quad
\mathcal{HF}_{Q(V_u)\{p_i,q_i\}}(\phi )(X,Z)\to
\mathcal {HF}^{\prime}_{Q(V_u)\{p_i,q_i\}}(\phi )(X,Z),
\end{equation}
where the $\phi (|X|,V)$ may depend also on 
$|X|$, or, $|X^{(a)}|$ and $|X^{(b)}|$. 
\end{theorem}
The proof is the same as for the constant basis case. Since both
basis', 
$\mathbf Q(V_u)$ and
$\mathbf Q^\prime (V_u)$,
are orthonormal, the proof of the
invariance of the radial functions, Euclidean Laplacian $\Delta_X$, and
projections $\Pi^{(n)}_X$ is even simpler as in the previous case.
Let it also be mentioned that versions (\ref{io4})-(\ref{io6}) can not
be introduced for the changing basis case because the 
$Q$'s must be in the intersection of all real
subspaces
$span_{\mathbb R}\mathbf Q(V_u)$.

Now we are ready to prove the first main theorem in this paper.
\begin{theorem}
Both in the constant and the changing basis cases, the
$\kappa_{\mathbf Q}$
intertwines the complete Laplacians $\Delta$ and $\Delta^\prime$.
\end{theorem}
\begin{proof}
There is proved above that both the Laplacian $\Delta_X$ and  
operators defined by multiplication with radial functions are
intertwined by the  
$\kappa_{\mathbf Q}$. The other parts of the Laplacians are 
also intertwined because of the following identities. 
\begin{eqnarray}
\label{MF}
\kappa_{\mathbf Q}:&  
\mathbf M\mathcal F_
{\mathbf Q(V_u)\{p_i,q_i\}}
(\phi )=
\mathcal F_
{\mathbf Q\{p_i,q_i\}}
((q-p)|V|\phi )
\\
&\to
\mathcal F^\prime_
{\mathbf Q\{p_i,q_i\}}
((q-p)|V|\phi )=
\mathbf M^\prime\mathcal F^\prime_
{\mathbf Q\{p_i,q_i\}}
(\phi ),
\nonumber
\\
\label{Delta_ZF}
\kappa_{\mathbf Q}:&  
\Delta_Z\mathcal F_
{\mathbf Q\{p_i,q_i\}}
(\phi )=
\mathcal F_
{\mathbf Q\{p_i,q_i\}}
(-|V|^2\phi )
\\
&\to
\mathcal F^\prime_
{\mathbf Q\{p_i,q_i\}}
(-|V|^2\phi )=
\Delta_Z\mathcal F^\prime_
{\mathbf Q\{p_i,q_i\}}
(\phi )
\nonumber
\end{eqnarray}
where $p=\sum p_i$ and $q=\sum q_i$. These formulas remain true if 
$\mathcal F$
is exchanged for
$\mathcal{HF}$.
\end{proof}

\subsection{Constructions 
on sphere$\times$ball-type domains.}

For introducing the intertwining operators 
on sphere$\times$ball-type domains, one can start with version 
(\ref{io3}), where it is defined for a constant basis in terms of 
homogeneous 
harmonic polynomials of the X-variable defined on the ambient space. 
It is pointed out there, that
the operator is well defined without using a basis, but now, for each $n$,
consider  functions of the form 
$
\Pi^{(n)}_X(
\prod_{i=1}^{k/2}z^{p_i}_{V_ui}\overline z^{q_i}_{V_ui})
=|X|^{2n}\sigma_{V_u\{p_i,q_i\}}
$
such that, for a fixed $V_u$, functions
$
\sigma_{V_u\{p_i,q_i\}}
$
form a basis among the corresponding spherical harmonics defined on the
unit sphere of the X-space. Note that the dimension of $n^{th}$-order 
harmonic polynomials is less than the dimension of $n^{th}$-order
homogeneous polynomials, thus, not all polynomials from the latter set
are subjugated to the projection.
Anyhow, such choices for such basis' exist.  
All those $V_u$'s satisfying this property
form an everywhere dense open subset of the unit Z-vectors.
The functions whose Fourier transforms are considered in an
arbitrary version
of the intertwining operator have unique expansions, 
$
\sum
_{\{p_i,q_i\}}
\phi
_{\{p_i,q_i\}}
(|X|,V)\sigma
_{V_u\{p_i,q_i\}}
$,
by these spherical harmonics
and the intertwining operator defined by this representation is 
the same as for the original representation. These spherical harmonics
can be pulled back from the unit X-sphere to the considered X-sphere by
the central projection $\pi: X\to X_u=X/|X|$.  
Then, for any function $c(V)$, 
function
$
c(Y)\pi^*(\sigma_{V_u\{p_i,q_i\}})
$ 
defined on the sphere$\times$ball-type manifold is the restriction exactly 
of those functions
$
\phi 
_{\{p_i,q_i\}}
(|X|,V)\sigma_{V_u\{p_i,q_i\}}
$ 
for which
$
\phi 
_{\{p_i,q_i\}}
(R_X ,V)=c(V)
$
holds. 
Thus we have
\begin{theorem} 
Both in the fixed and changing basis cases, intertwining operator 
$\kappa_{\mathbf Q }$ induces a well defined action on  
functions defined by restrictions from the ambient 
manifolds onto the sphere$\times$ball submanifolds. 
This induced  
operator, 
$\tilde\kappa_{\mathbf Q }$,
is well defined
for all versions  
and is the same as the operator constructed
by a basis of the space of spherical harmonics.

This induced operator intertwines the Laplacians defined
on the sphere$\times$\-ball-type submanifolds. 
\end{theorem}
The proof of well-definedness for 
a basis of the space of spherical harmonics
is the same as on the ambient manifold. The Laplacian on the submanifold
differs from (\ref{Delta}) just by the terms $\Delta_X$ and $|X|$ 
which should 
be exchanged for $\Delta_{S_X}$ (which is the Laplacian on the X-sphere)
and $R_X$ (which is the radius of the X-sphere), respectively. Note that
the X-directional derivatives included into $\mathbf M$ concern directions
tangent to the sphere. Thus this term is the same as for the ambient space.
In other words, in order to have the Laplacian on the submanifold, just 
the radial Laplacian, $\Delta_r$, of the X-space
should be dropped from (\ref{Delta}). Since both the radial and
the complete Laplacians are 
invariant under the action of the ambient intertwining operator, also
the Laplacian on the submanifold is invariant.

\section{Domain and range of $\kappa_{\mathbf Q}$.}
The domain and range of the intertwining operators is determined by  
a function transformation which is noteworthy 
also without this application. 

\subsection{The dual Radon transform.} 

This transform was first investigated in \cite{sz6}, 
pages 264-266, where it is called {\it boomerang transform}. 
The results provided there include also an inversion formula, which,
by a new proof, was reestablished by 
\'A. Kurusa \cite{ku1, ku2}. 
He called the operator itself {\it dual Radon transform} which 
name better describes the area this transform belongs to.
We adopt this name, however, 
the following review proceeds with the author's original ideas.

Let 
$g_\theta (r)$
be a half-line parameterized by arc-length $r$ which has its endpoint, 
corresponding
to $r=0$, at the origin
$O$ of $\mathbb R^l$ and which is pointing to the point
$\theta\in S^{l-1}_0(1)$
of the unit sphere 
$S^{l-1}_0(1)\subset\mathbb R^l$
around the origin $O$. Then 
$(\theta ,r)$ 
serve as polar coordinates for the points, $Z$, of 
$\mathbb R^l$. These denotations indicate that this transform will be
used on the Z-space of H-type groups. If 
$f(\theta ,r)$
is a continuous function defined on  
$\mathbb R^l$,
then, for each fixed 
$\theta_0$,
it determines a cylindrical function
$f^c_{\theta_0}(Z)$
defined on the unique half-space whose perpendicular 
projection onto the line  spanned by 
$g_{\theta_0}(r)$
is equal to this half-line. If the projection of $Z$ is the point having
the polar coordinates
$(\theta_0 ,r)$,
then, by definition, 
$f^c_{\theta_0}(Z)=f(\theta_0 ,r)$. By considering this construction
for each 
$\theta$,
one can define the function-valued function 
$\theta\to f^c_{\theta}$.
The dual Radon transform, $f\to f_\tau$, is defined by the integral
\begin{equation}
\label{tau}
f_\tau :=\int_{S^{l-1}}
f^c_{\theta}d\theta ,
\end{equation}
which can be written also in the form
\begin{equation}
f_\tau (Z):=\int_{
\langle\theta ,Z\rangle
\geq 0}
f (\theta ,
\langle\theta ,Z\rangle )
d\theta .
\end{equation}
Apply Thales' theorem to the last formula to see that the transform
is defined by the integral of $f$   
on the sphere of
diameter $[0,Z]$  
by the measure $d\theta$. 
Note that this 
measure differs from the canonical measure of the Thales sphere.
By (\ref{tau}) and Fubini's theorem we have:
\begin{lemma} Let $f(Z)$ be an arbitrary continuous function and
let $\mu$ be a continuous function with compact support in $\mathbb R^l$.
Then the integral formula
\begin{equation}
\label{rdual}
\int_{\mathbb R^l}
f_\tau (Z)\mu (Z)dZ=
\int_{S^{l-1}}
\int_0^\infty
f(\theta ,r)
\mu_R(\theta ,r)
drd\theta
\end{equation}
holds, where
$
\mu_R(\theta ,r)
$
is the Radon transform of $\mu$ defined by the integrals of this
function on the hyperspaces intersecting 
$g_\theta$ at the points having the polar coordinates 
$
(\theta ,r)
$
perpendicularly.
\end{lemma}
Formula (\ref{rdual}) reveals that the considered transform is dual
to the Radon transform, indeed. The main result in this section is:
\begin{theorem} 
Let 
$f_\tau(\theta ,r)$
be an arbitrary function of class $C^{2m}$ with compact support in 
$\mathbb R^l$, where $l=2m+1$ is odd or $l=2m$ is even. Then $f_\tau$
has an inverse, $f$, regarding the dual Radon transform, 
which is of the form
\begin{equation}
\label{inv1}
f={(-1)^m\over (2\pi )^{2m}} 
( (f_\tau )_R)^{(2m)},
\quad {\rm if}\quad l=2m+1,
\end{equation}
where $(2m)$ means the $2m$th derivative of the functions with respect
to $r$, resp.
\begin{equation}
\label{inv2}
f={(-1)^m(l-1)!\over (2\pi )^{2m}} ( (f_\tau )_R)^{[2m]},
\quad {\rm if}\quad l=2m,
\end{equation}
where 
$\varphi^{[2m]}(r)$ 
is defined for a function $\varphi (t)$ on
$\mathbb R$ by
\begin{eqnarray}
\varphi^{[2m]}(r):= \int_0^\infty
{1\over t^{2m}}\big(
\varphi (r+t)+ 
\varphi (r-t) 
\\
 -2\big[
\varphi (r)
+{t^{2}\over 2!} 
\varphi^{\prime\prime} (r)
+\dots  
+{t^{2m-2}\over (2m-2)!} 
\varphi^{(2m-2)} (r)
\big]\big)dt.
\nonumber 
\end{eqnarray}
\end{theorem}
\begin{proof}
Let $\mu$ be a function of class $C^{2m}$
with compact support in $\mathbb R^{2m+1}$. 
From (\ref{rdual}) we get
\begin{eqnarray}
\int_{\mathbb R^l}[( (f_\tau 
)_R)^{(2m)}]_\tau (Z)\mu (Z)dZ=\\
\int_{S^{l-1}}
\int_0^\infty
{(-1)^m\over (2\pi )^{2m}} 
( (f_\tau )_R)^{(2m)}
(\theta ,r)
\mu_R
(\theta ,r)
drd\theta =
\\
\int_{S^{l-1}}
\int_0^\infty
 (f_\tau )_R
(\theta ,r)
\big(
{(-1)^m\over (2\pi )^{2m}} 
(\mu )^{(2m)}_R (\theta ,r)\big)
drd\theta=\\
\int_{\mathbb R^l}f_\tau (Z)
\big(
{(-1)^m\over (2\pi )^{2m}} 
(\mu )^{(2m)}_R\big)_\tau (Z)dZ=
\int_{\mathbb R^l}f_\tau (Z)\mu (Z)dZ.
\nonumber
\end{eqnarray}
In the last step the well known Radon inverse formula is used. Since $\mu$
is arbitrary, formula 
(\ref{inv1})
is established. Formula
(\ref{inv2})
can be established in the same way.
\end{proof}
Let it be mentioned that all non-trivial isospectrality examples 
constructed in this paper arise from odd dimensional Z-spaces. 
Thus, only formula (\ref{inv1}) applies to these cases.

\subsection{The domain of the intertwining operators.} 

The above theorem is used to prove that the function space generated by
functions of the form
$
\mathcal F_{\mathbf Q\{p_iq_i\}}(\phi )(X,Z)
$
contains all functions 
$P(X)f(Z)$, 
where $P(X)$ is a complex
valued polynomial and $f(Z)$ is a smooth function of compact support on
the Z-space. Thus, by using appropriate limiting procedures,
the whole standard $L_{\mathbb C}^2$-Hilbert space on the 
$(X,Z)$-space can be generated in this way. Note that the Z-Fourier 
transform of $\phi$ is ``twisted'' with the polynomials appearing 
in the formula which depend, beside $X$, 
also on $V_u$. By this reason, it is called
also {\it twisted Z-Fourier transform}. It is this feature what makes
the constructions of the above functions highly non-trivial.

The proof of the above statement needs some preparations. For a fixed unit 
Z-vector 
$V^0_u$ 
and positive number $\delta$, let 
$T_\delta (V^0_u)$ 
be the tube of radius $\delta$ around the half-line 
$g_{V^0_u}(r)$. By the standard definition, 
it is the union of those discs, 
$D^{(l-2)}_\delta(r)$,
of radius $\delta$ about the points of the half-line which 
intersect the half-line perpendicularly. The characteristic functions 
of this tube and the half-line, 
defined on the whole Z-space resp. line determined by
$g_{V^0_u}(r)$, 
are denoted by   
$\chi_{\delta V^0_u}(Z)$
and
$\chi_{g_{V^0_u}}(t)$ 
respectively. 
Then, 
\begin{eqnarray}
\lim_{\delta\to 0}{1\over 
Vol(D^{(l-2)}_\delta )}
\mathcal F_{\mathbf Q\{p_iq_i\}}(
\chi_{\delta V^0_u}
\phi )(X,Z)
=\prod_{i=1}^{k/2}z^{\prime p_i}_{V^0_ui}(X)
\overline z^{\prime q_i}_{V^0_ui}(X)\mathbb 
L_{V_u}^c
(\phi_{V^0_u})(Z),
\nonumber 
\end{eqnarray} 
where function 
$
L_{V_u}(\phi_{V^0_u})(t)=
\chi_{g_{V^0_u}}(t)
Fou_{\pm V_u}(
\chi_{g_{V^0_u}}
\phi_{V^0_u})(t),
$ 
defined on the whole line spanned by $V_u^0$ 
and parameterized by $t$
satisfying $t(V_u)=1$, vanishes for $t < 0$ 
and it is the Laplace transform of
$
\phi_{V^0_u}(r)
$
defined on the half-line 
$g_{V^0_u}(r)$.
In the latter formula, this function is described in terms of the
1-dimensional Fourier transform
$
Fou_{\pm V^0_u}
$ 
defined on the whole line spanned by $V_u^0$. 

Note that the above function does appear as a product of 
X- and Z-depending functions. In the next averaging process they are 
used to 
construct 
$P(X)f(Z)$ such that one considers the same polynomial $P(X)$ for all
$V_u$ and, in the end, the $f$ appears as the dual Radon transform of
an appropriate function defined on the Z-space.
 
For a complex valued polynomial, $P(X)$, let 
$
P_{V^0_u}(X) 
$
be a representation of the polynomial in terms of the complex structure
$J_{V^0_u}$. Furthermore, let 
$
P_{V^0_u}(X,V_u) 
$
be the function defined by replacing
$V^0_u$ with general $V_u$ in  
$
P_{V^0_u}(X). 
$
Therefore,
$
P_{V^0_u}(X)= 
P_{V^0_u}(X,V^0_u) 
$
holds, but for other $V_u$'s there are other polynomials defined. 
Denotation 
$\mathcal F_{\mathbf Q
P_{V^0_u} 
}$
means that
$\prod_{i=1}^{k/2}z^{\prime p_i}_{V^0_ui}(X)
\overline z^{\prime q_i}_{V^0_ui}(X)$
is replaced by 
$
P_{V^0_u}(X,V_u) 
$
in the above formulas. Then we have:
\begin{eqnarray}
\lim_{\delta\to 0}
\int_{S^{l-1}}
{1\over Vol(D^{(l-2)}_\delta )}
\mathcal F_{\mathbf QP_{V^0_u}}
(\chi_{\delta V^0_u}
\phi )(X,Z)dV^0_u=\\
P(X)
\int_{S^{l-1}}
\mathbb L^c_{V^0_u}(\phi_{V^0_u})(Z)dV_u^0=
P(X)
(\mathbb L_{\pm V_u}(\phi ))_\tau (Z),
\nonumber
\end{eqnarray} 
where
$
\mathbb L_{\pm V_u}(\phi ) 
$
denotes the function
defined on the whole line spanned by $\pm V_u$ 
by the Laplace transforms of functions
$\phi_{V_u}(r)$
resp. 
$\phi_{-V_u}(r)$.
Thus we have:
\begin{theorem}
For given polynomial $P(X)$ and smooth function $f(Z)$ of compact support
the product $P(X)f(Z)$ is limit of convergent sequences of functions
belonging to the domain of an intertwining operator. This sequence is 
constructed by the above method, where function
$\phi (V) =\mathbb L^{-1}_{\pm V_u}(f_{\tau^{-1}})(V)$
is derived from $\phi$ by the inverse dual-Radon resp. 1-dimensional
Laplace transforms defined above on the corresponding half-lines. 
(The inverse formula for the Laplace transform is called Mellin's formula.
An alternative version is the so called Post's formula.)  

The same proof yield those versions of the theorem when function 
$f$ is of the form 
$f(|X|,Z)$
such that, for any fixed $|X_0|$, function
$f(|X_0|,Z)$
is of compact support on the Z-space, or, 
when this problem is considered on
a sphere$\times$ball-type domain and 
$P(X)$ is replaced by its restriction,
$\tilde P(X)$, onto the sphere and $f(Z)$ is 
the same function as before.
In these cases 
$P(X)f(|X|,Z)$
resp.   
$\tilde P(X)f(Z)$
are in the domain of the intertwining operator. 
In both cases, functions
$\phi (|X|,V)$ resp. 
$\phi (V)$
can be found by the same inverse operations defined on the Z-space.
\end{theorem}

\section{Intertwining of the boundary conditions.}
The most important tool applied in establishing 
the intertwining of the boundary conditions is a theory developed
for one- and two-pole functions.
\subsection{Formulas for one- and two-pole functions.}\label{tdc}

For a unit Z-vector $Z_0$ and $Q\in\mathbb R^k$, denotation
$X_{QZ_0}$ means that this X-vector is in the subspace spanned
by $Q$ and $J_{Z_0}(Q)$. On the plane 
$P(Q,J_{Z_0}(Q))$ 
spanned by these two vectors
the polar coordinates $(|X_{QZ_0}|,\alpha )$ are defined such that
$
\alpha (Q)=0,
\alpha (J_{Z_0}(Q))=\pi /2
$ hold. 
By the restriction 
$\alpha\leq\pi$ 
imposed for all $Z_0$, 
one has a spherical
coordinate system on the $(l+1)$-dimensional space, $S_Q$, spanned by 
$Q$ and all  
$J_{Z_0}(Q)$. Thus these $\alpha$ parameter lines are half circles running
in the half-plane
$
P^+(Q,J_{Z_0}(Q))\subset
P(Q,J_{Z_0}(Q))
$ 
bounded by $\mathbb RQ$ and containing
$
J_{Z_0}(Q)
$.

Function $\Theta_Q$ can be described by this coordinate 
system as follows. If the orthogonal projection, $X_Q$, of $X$
onto $S_Q$ is in $P^+(Q,J_{Z_0}(Q))$, then
\begin{equation}
\label{theta}
\Theta_Q(X,V_u)=\langle Q,X_Q\rangle +\mathbf i\langle[Q,X_Q],V_u\rangle
=|X_Q|(\cos\alpha +\mathbf i \langle Z_0,V_u\rangle\sin\alpha ).
\end{equation}
Powering performed in $\Theta_Q^p\overline\Theta_Q^q$ yield:
\begin{theorem}
On those vectors, $X$, whose projections onto $S_Q$ 
fall onto a fixed $\alpha$-half-circle around
the origin of the half-plane 
$P^+(Q,J_{Z_0}(Q))$, a 1-pole function
$
\mathcal F_{Qpq}(\phi )
$
has the form: 
\begin{eqnarray}
\mathcal F_{Qpq}(\phi )(X,Z)=
\int_{\mathbb R^l}
\phi (|X|,V)
\Theta_{Q}^p(X,V_u)\overline\Theta^q_{Q}(X,V_u)
e^{\mathbf i\langle Z,V\rangle}
dV=
\label{eqn1}
\end{eqnarray}
\begin{eqnarray}
\sum_{s=0}^{p+q}|X_Q|^{p+q}
\cos^{p+q-s}\alpha \sin^s\alpha 
\int_{\mathbb R^l}
A_{spq}
\langle Z_0,V_u\rangle^s
\phi (|X|,V)
e^{\mathbf i\langle Z,V\rangle}
dV=\label{eqn2}
\end{eqnarray}
\begin{eqnarray}
\sum_{s=0}^{p+q}|X_Q|^{p+q}
\cos^{p+q-s}\alpha \sin^s\alpha 
(A_{spq}
(-\mathbf i)^s
\partial^s_{Z_0})\int_{\mathbb R^l}
\phi (|X|,V)|V|^{-s}
e^{\mathbf i\langle Z,V\rangle}
dV.
\label{eqn3}
\end{eqnarray}
Function $\phi_s(X,V)=\phi (|X|,V)|V|^{-s}$, 
whose Z-Fourier transform, $\tilde\phi_s$, appears 
as the last integral term of (\ref{eqn2}), 
is derived from $\phi$ such that it  
depends just on $|X|$ and the Z-variable. Term behind 
$\sin^s\alpha $ is denoted by
$\tilde A_{spq}(
|X|,Z_0,Z)$.
If both
$|X|$ and $Z$ are fixed, then
$
\tilde A_{spq}
$
is constant for those $X$'s which project onto the half circle determined
for $Z_0$ by the parameter-range 
$\alpha\leq\pi$. 
On $S_Q$, whose points are denoted by $X_Q$, 
this function appears in the form
\begin{equation}
\label{F/A}
\sum_{s=0}^{p+q}|X_Q|^{p+q}
\cos^{p+q-s}\alpha \sin^s\alpha \,
\tilde A_{spq}(|X_Q|,Z_0,Z),
\end{equation} 
where, for fixed values of
$|X|$ and $Z$, the
$
\tilde A_{spq}
$
is an $s^{th}$-order polynomial which can be described in terms of the 
unit vectors $J_{Z_0}(Q_u)$, where $Q_u=Q/|Q|$, as follows.
 
Originally, this polynomial can explicitly be determined 
on the Z-space by the expansion 
$
\langle Z_0,V_u\rangle^s
=\sum_{j=0}^sB_j\sigma_{Z_0}^{s-j}(V_u)$
in terms of the spherical harmonics
$\sigma_{Z_0}^{s-j}(V_u)$. 
Since, for any fixed $Z_0$, function
$
\langle Z_0,V_u\rangle^s
$
defined on the unit Z-sphere is radial
about the center $Z_0$, thus also the spherical harmonics are radial
about $Z_0$ and the convolutions with them are nothing but the
projections onto the corresponding subspaces of spherical harmonics. Thus,
\[
\tilde A_{spq}
(|X|,Z_0,Z)
=\sum_jB_j
\varphi^{(s-j)}
(|X|,Z_0,Z),
\]
where
$
\varphi^{(s-j)}
$
is the corresponding spherical harmonics appearing in the expansion of
$
\mathbf\varphi (|X|,Z,V_u)=\int_0^\infty A_{spq}\phi_s(|X|,V_u,r)
e^{r\mathbf i\langle Z,V_u\rangle}dr
$
which function is defined, for fixed $Z$ and $|X|$, by integrals 
with respect to $dr$ defined for the polar coordinate system $(V_u,r)$. 
 
But this function depends on $X$. In its final form, it can be viewed 
such that the function determined on the Z-space defines, first, a
0-homogen\-eous function
on the equator plane 
$E_{Q_u}$ 
spanned by the X-vectors $J_{Z_0}(Q_u)$. Then, it extends onto $S_Q$
such that, on an $X_Q$, it takes the value determined by $Z_0$ 
if and only if
$
X_Q\in
P^+(Q,J_{Z_0}(Q))
$.

For other points, which are outside of $S_Q$, the function is determined
by this function and projection onto $S_Q$. Note that the
$
\tilde A_{spq}
$
is defined for $X$ and not for $X_Q$, meaning that, 
instead of $|X_Q|$, 
function $\phi$
involves 
$|X|$ 
to this term.
The latter parameters stand in front of the formula
and are in connection with the trigonometric polynomials. 
\end{theorem}

Such formulas can be established also for
$
\mathcal {HF}_{Qpq}(\phi )
$. 
Functions $\overline\Theta_Q^q$ resp. $\Theta_Q^p$ are homogeneous 
harmonic polynomials of the X-variable, thus in cases satisfying
$p=0$ or $q=0$, the function in (\ref{eqn1}) is nothing but 
$
\mathcal {HF}_{Qpq}(\phi )
$. 
If $pq\not =0$, there are
new terms,
$|X|^{2r}\Theta_Q^{p-r}\overline\Theta_Q^{q-r}$,
appearing in the X-harmonic polynomial 
$\Pi_X(\Theta_Q^{p}\overline\Theta_Q^{q})$. 
For each $r$, an additional sum
shows up both in (\ref{eqn2}) and (\ref{eqn3}). Comparing the first 
and the $r^{th}$ sums, the $p+q$ is exchanged for 
$p+q-2r$, which is the greatest
possible value for $s_r$ in the sum. Such a new sum can be combined
with the first one, where $r=0$, by multiplying the $r^{th}$ sum by
$1=(\cos^2\alpha +\sin^2\alpha )^r$. Thus one gets trigonometric
polynomials appearing in the first sum. 
By collecting the terms belonging to the same trigonometric 
polynomial, the first term behind the integral sign of (\ref{eqn2}) 
is exchanged for a polynomial of the form
$
P_{spq}(|X|,\langle Z_0,V_u\rangle )=
\sum_{r=0}^sA^{(s-r)}_{spq}|X|^r
\langle Z_0,V_u\rangle^{s-r}
$,
resulting the integral terms

\begin{eqnarray}
\label{P}
\tilde P_{spq}(|X|,Z_0,Z )=
\int_{\mathbb R^l}
P_{spq}
(|X|,Z_0,V_u)
\phi (|X|,V)
e^{\mathbf i\langle Z,V\rangle}
dV,
\\
\tilde P_{spq}(|X|,\partial^{s-r}_{Z_0},\tilde\phi_{s-r})=  
\sum_{r=0}^sA^{(s-r)}_{spq}(-\mathbf i)^{s-r}|X|^r
\partial_{Z_0}^{s-r}\tilde \phi_{s-r}
\nonumber 
\end{eqnarray}
behind $\sin^s\alpha$
in formulas (\ref{eqn2}) resp. (\ref{eqn3}). Note that constants 
$A^{(s-r)}_{spq}$
are built up by but not equal to the constants
$A_{spq}$. Thus we have:
\begin{theorem}
On $S_Q$, function
$
\mathcal {HF}_{Qpq}(\phi )
$ 
appears in the form
\begin{equation}
\label{F/P}
\sum_{s=0}^{p+q}|X_Q|^{p+q}
\cos^{p+q-s}\alpha \sin^s\alpha \,
\tilde P_{spq}(|X_Q|,Z_0,Z),
\end{equation} 
where 
$
P_{spq}
$
is explicitly described in (\ref{P}). 
For fixed values of
$|X|$ and $Z$ also this term
is an X-depending $s^{th}$-order polynomial which appears in the same
form as 
$
\tilde A_{spq}
$
does. But this one has also lower order terms,
$
\langle Z_0,V_u\rangle^{s-r},
$
beneath the main term. For other points not being on $S_Q$ 
also this function is determined
by projections onto $S_Q$, in which case it
is defined in terms of $|X|$ and not $|X_Q|$.
\end{theorem}

The above constructions restricted
onto spheres $S_{R_X}$ provide the formulas 
on sphere$\times$ball-type domains. In this case function 
$|X|$ is constant, thus
functions $\tilde\phi_{s-r}$ depend just on $Z$ and $Z_0$. Let it be 
pointed out again that the latter 
variable is involved by the assumption
$X_Q\in P^+(Q,J_{Z_0}(Q))$, 
i. e., it is determined by $X$ over which
the Fourier transform in the Z-space is performed. In other words, it is
an X-depending function whose precise denotation would be $Z_0(X)$. 

Since the Z-balls, $B_{R_Z}(X)$, where 
$X\in S_{R_X}$, 
are naturally 
identified on this trivial ball-bundle, 
they determine the same functions in the Z-space for all those $X$'s which
project onto the half-plane 
$P^+(Q,J_{Z_0}(Q))$.
Thus functions 
$
\mathcal F_{Qpq}(\phi )
$
resp.
$
\mathcal {HF}_{Qpq}(\phi )
$
appear in the form 
(\ref{F/A}) resp. 
(\ref{F/P}) such that $|X|=R_X$ is constant in this case.

Later on, we need these functions described also on circles, 
$
C=P_2\cap
S_{R_X}
$,
which are represented as intersections of 2-dimensional linear subspaces,
$P_2$, with $S_{R_X}$. If  
$
Q\in
P_2
$,
these functions are perfectly described by the above formulas also on these
circles. Therefore, we suppose
$
Q\not\in
P_2
$. First also suppose that
$
P_2\subset
S_{Q}
$
holds, in which case the computations below are carried out on the
3-space,
$
S_{Q3}
$, 
spanned by
$
P_2
$
and
$
\mathbb RQ
$.
This space intersects 
$S_{R_X}$ 
at the 2-sphere denoted by
$S_{R_X2}$. 
The north-pole, $O$, of this sphere is cut out by the ray
$\mathbb R_+Q$. The north-pole, $O_C$, on the circle is defined by 
the closest
point to $O$. Let $\alpha_C$ be the angle parameterization of $C$ with 
origin $O_C$ such that on both sides  
$0\leq \alpha_C\leq\pi$ 
hold. The angle between $C$ and the 
great circle $CE_O$ with center $O$ (equator) 
on the 2-sphere   
$S_{R_X}2$ 
is denoted by 
$\beta_C$. It is uniquely determined by the assumption 
$0\leq \beta_C\leq\pi/2$. 
For a point $P\in C$ satisfying 
$0\leq \alpha_C\leq\pi/2$, let $\tilde C_P$ be the great circle connecting
$O$ and $P$, which intersects $CE_O$ at a point $N$ perpendicularly. If
$M_C=C\cap CE_O$, then the spherical sine theorem applied to the 
right spherical triangle $PNM_C$ yields
\begin{equation}
{\sin({\pi\over 2}-\alpha (P))\over   
\sin \beta_C}=   
{\sin({\pi\over 2}-\alpha_C (P))\over   
\sin {\pi\over 2}}\,\,
\Rightarrow   \,\,
\cos \alpha (P)=   
\sin \beta_C   
\cos \alpha_C (P)
\end{equation}
This equation along with 
$
\sin \alpha (P)=\sqrt{1-      
\sin^2 \beta_C   
\cos^2 \alpha_C (P)}
$
imply that functions 
(\ref{F/A}) 
and
(\ref{F/P})
restricted onto $C$ can be expressed in terms of
$
\cos\alpha_C (P)
$.

It is a very important issue to understand the precise appearance of 
functions
$
\tilde A_{spq}
$
and
$
\tilde P_{spq}
$
on these circles. They appear as polynomials on the equator circle 
$CE_O$ but on $C$ they appear as functions which are
pulled back from the equator to the $C$ by the 
central projectivity 
$\tau_O:C\to CE_O$. 
This $\tau_O$ can be explicitly computed
as follows. 

Parameterize both $P_2$ and $E_O$
by complex numbers $z$ and $z^\prime$ respectively such that
these coordinate systems have common imaginary axis 
$\mathbb R_+\mathbf i\subset P_2\cup E_O$
and $z=1$ and $z^\prime=1$
correspond to $P=O_C$ and $N_P=\tilde C_P\cap CE_O$ respectively, where
$\tilde C_P$ is the great circle connecting
$O$ and $P$. This circle intersects $CE_O$ at $N_P$ perpendicularly.
Pick up also such unit complex numbers, $u$ and $u^\prime$, between
the units and the imaginary numbers which are corresponded to each other
by the $\tau_O$. 

Actually, the 
$\tau_O$
is a real projectivity between the two projective lines
$\tilde C$ and $\tilde{CE}_O$ 
defined by identifying the antipodal points on the great circles
$C$ and ${CE}_O$. Thus, it can be described in terms of the real cross
ratio defined on these projective lines. But this real one is the same
as the complex cross ratio defined on the complex planes if the points
are laying on the same half-circle. (This statement is well known in
conform geometry of 2-spheres.) Also note that the common imaginary
numbers cut the great circles into half-circles which are corresponded
to each other by the
$\tau_O$,
therefore, this projectivity can be described in terms of
the complex cross ratio by the relation
$
(z^\prime =\tau_O(z),1^\prime,\mathbf i^\prime,u^\prime)
=(z,1,\mathbf i,u). 
$ 
This equation describes the 
$\tau_O$  
as a fractional linear function (M\"obius transform) of the form 
$z^\prime (z)=(az+b)/(cz+d)$. 
Since such a function preserves the circles and the three
corresponding points
$(1^\prime,\mathbf i^\prime,u^\prime)$ and $(1,\mathbf i,u)$
are on $CE_O$ and $C$, the transformation is really defined between the
two circles. Thus we have:
\begin{lemma} The
$\tau_O$
is a fractional linear function (M\"obius transform) between the two
great circle which pulls back a trigonometric function 
$\cos_u\alpha=(1/2)(uz+\overline u\overline z)$,
defined on the complex plane $E_0$ by a fixed complex unit $u$, 
to the trigonometric rational function  
$(a\cos_u+b)/(c\cos_u+d)$
defined on $C$. All trigonometric polynomials can be generated on $CE_0$
by the functions $\cos_u$, therefore, all trigonometric polynomials
are pulled back to a trigonometric rational function defined on the other
great circle. Functions 
$
\tilde A_{spq}
$
and
$
\tilde P_{spq}
$
can be described by the trigonometric functions 
$
\varphi^{(s-j)}
$
defined on the unit sphere of $E_0$. Thus they are trigonometric 
polynomials on $CE_0$ which pull back to trigonometric rational
functions defined on $C$.
\end{lemma}

For a $P_2$ which is not subspace of $S_Q$ these functions can be 
determined by projecting it into $S_Q$. Almost every plane projects to a
plane of $S_Q$.  If 
$P_2^\prime$ 
is the projected
plane, then the sought functions on $P_2$ are the pull-back's of the
corresponding functions defined for 
$P_2^\prime$.
Note that $\cos\alpha_C^\prime$
arises from a linear function, therefore, so does the pull-back 
function whose kernel is the pull-back of the line (linear subspace) 
connecting the points 
$M_{C^\prime}$
and
$-M_{C^\prime}$
defined above for $C^\prime$. Then the closer midpoint, $O_C$, to the $O$, 
which is between the pull-back-points
$M_{C}$
and
$-M_{C}$ on $C$,  
is called the north-pole on $C$. 
This point determines the parameterization $\alpha_C$. The functions
restricted onto $C$ are described in terms of this parameter.
Thus we have
\begin{theorem}
\label{restric}
When the constructions are restricted
onto a fixed sphere $S_{R_X}$, functions 
$
\mathcal F_{Qpq}(\phi )
$
resp.
$
\mathcal {HF}_{Qpq}(\phi )
$
appear in the form 
(\ref{F/A}) resp. 
(\ref{F/P}) such that functions
$
\tilde A_{spq}
$
and
$
\tilde P_{spq}
$
involve the constant $|X|=R_X$.
In the following statement denotation
$
\tilde R_{spq}
$
can be replaced by any of these two functions.

On a circle, 
$
C=P_2\cap
S_{R_X}
$,
represented by intersection of a 2-dimensional linear subspace
$P_2$ with $S_Q$, these functions are of the form  
\begin{equation}
\label{circ}
\sum_{s=0}^{p+q}K_C^{p+q}
(\sin\beta_{C^\prime} 
\cos
\alpha_C) 
^{p+q-s}
(1-\sin^2 
\beta_{C^\prime
}   
\cos^2 \alpha_C)^{s\over 2}
\tilde R_{spq}(Z_0,Z),
\end{equation} 
where
$
C^\prime
$
is the projected circle cut out by the projected 2-space  
$
P_2^\prime
$
which intersects the equator $E_O^\prime$ on the projected space
$
S_{Q3^\prime}
$
at angle
$
\beta_{C^\prime}   
$. Constant $K_C^{p+q}$ is due to the fact that the pulled back linear
functions are restricted to a circle in this process. 
For circles in $S_Q$ this constant is $K_C=R_X$.

For any fixed $Z$ and $|X|$, function
$
\tilde R_{spq}(Z_0,Z)
$
defines a trigonometric rational function on $C$ which is the pull back
of an $s^{th}$-order polynomial with such a combined map, where the first
map, $\tau^\prime_O$, is a M\"obius transform between 
$C^\prime$ to $CE_O^\prime$ and the second one takes $C$ onto $C^\prime$
by a projection.  
\end{theorem}

We need these theorems in the more general case when, instead of
$Q$, one considers a pair, 
$(Q^{(a)},Q^{(b)})$, 
of vectors and functions 
$\phi (|X|,V)$,\,
$
\Theta^{p}_{QV_u}\overline \Theta^{q}_{QV_u}
$,\,
$
\Pi^{(p+q)}_X(..)
$
are exchanged for the following ones 
\[
\phi (|X^{(a)}|,|X^{(b)}|,V),\,\,
\Theta^{p_a}_{Q^{(a)}V_u}\overline \Theta^{q_a}_{Q^{(a)}V_u} 
\Theta^{p_b}_{Q^{(b)}V_u}\overline \Theta^{q_b}_{Q^{(b)}V_u},\,\,
\Pi^{(p_a+q_a)}_{X^{(a)}}(..) 
\Pi^{(p_b+q_b)}_{X^{(b)}}(..),
\]
respectively.
Such functions are called 2-pole functions which 
can be investigated in two ways. They can be considered
either on subsets
$(X_F^{(a)},
X^{(b)})$
defined by a fixed 
$X_F^{(a)}$,
or, on the similarly defined subsets
$(X^{(a)},
X_F^{(b)})$. Because of the exact similarities, only the first case
should be described, when, functions
\begin{equation}
\label{X^aV}
\phi (|X_F^{(a)}|,
|X^{(b)}|,V),\,\,
\Theta^{p_a}_{Q^{(a)}V_u}\overline \Theta^{q_a}_{Q^{(a)}V_u},\,\,
\Pi^{(p_a+q_a)}_{X^{(a)}}(..) 
\end{equation}
depend (non-trivially) just on
$
|X^{(b)}|
$
and
$
V
$.
On a circle
$(X_F^{(a)},
C^{(b)})$
the considered functions appear in the form (\ref{circ}) where the last
function is defined by those listed in (\ref{X^aV}) and the other
functions are defined on
$
\mathbf v^{(b)}
$.

\subsection{Intertwining of the Dirichlet conditions.}
The Dirichlet Intertwining Theorem will be established, first, 
for a constant basis, $\mathbf Q_F$. The changing basis
case will be traced back to this first one.
 
Observe that functions 
$ 
\cos^{p+q-s}(\alpha )\sin^s(\alpha )
$,
satisfying
$0\leq s\leq p+q$
are linearly independent, furthermore, for any fixed $Z_0$ and $Z$, 
function  
$\tilde A_{spq}$ is constant on the $\alpha$-parameter line determined
by $Z_0$. These two statements yield the following theorem obviously.
\begin{theorem}
A function 
$\mathcal{F}_{Qpq}(\phi )$ satisfies the Dirichlet condition 
at the boundary points $(X,Z)$  
if and only if functions
$\tilde A_{spq}$ 
vanish on the sphere $S_{R_Z}$, for all $Z_0(X)$ and 
$0\leq s\leq p+q$. 
Regarding
$\mathcal{HF}_{Qpq}(\phi )$, this condition is 
$\tilde P_{spq}=0$, for all
$0\leq s\leq p+q$ and $Z_0(X)$ at any boundary point $Z\in S_{R_Z}$.

For fixed $Q$ and natural numbers $p$ and $q$, function spaces
$
\mathbf \Phi_{Qpq}
$
resp. 
$
\mathbf \Xi_{Qpq}
$
are defined by the $L^2$ function spaces spanned by functions of the form
$\mathcal F_{Qpq}(\phi )(X,Z)$
resp. $\mathcal{HF}_{Qpq}(\phi )(X,Z)$,
where $\phi(|X|,V)$ (which depends, non-trivially, just on $V$ on 
sphere$\times$sphere-type manifolds)
can be an arbitrary $L^2$-function. For fixed $Q$ but running $p$ and $q$, 
all these spaces sum up to the total space
$
\mathbf \Phi_Q
=\sum_{p,q}
\mathbf \Phi_{Qpq}=
\mathbf \Xi_Q=\sum_{p,q}
\mathbf \Xi_{Qpq}
$.
Then, for functions $\varphi_Q$ resp. 
$\varphi_Q^\prime =\kappa_{\mathbf Q_F}(\varphi_Q)$ from
$
\mathbf \Phi_Q
$
resp.
$
\mathbf \Phi^\prime_Q
$
the Dirichlet condition is satisfied always simultaneously. Actually,
the intertwining operator
$\kappa_{\mathbf Q_F}:
\mathbf \Phi_Q \to
\mathbf \Phi^\prime_Q
$
between these total spaces is induced 
by a point transformation of the form
$(T_Q(X),id_Z)$, where the $T_Q$ is 
an orthogonal transformation on the X-space, depending on $Q$. 
\end{theorem}

Such simple proof can be given only for total 
spaces defined by a fixed pole
$Q$. The proof is much more difficult on the complete 
$L^2$-Hilbert space, which can be represented both as  
$
\sum_{Q\in span_{\mathbb R}(\mathbf Q_F)}
\mathbf \Phi_Q, 
$
and
$
\sum_{Q\in span_{\mathbb R}(\mathbf Q_F)}
\mathbf \Phi^\prime_Q,
$
i. e., by the sums of all functions defined by all poles, $Q$, which 
are in the real span of $\mathbf Q_F$. The main idea of such an extension
is as follows.

Suppose that a function
$\varphi (X,Z)=\sum_{Q_i}\varphi_{Q_i}(X,Z)$
satisfies the Dirichlet condition.
Decompose each $\Theta_Q$ in the form 
$\Theta_Q(X,V_u)=
\langle Q,X\rangle+\mathbf i
\langle J_{V_u}(Q),X\rangle
$ and $\overline\Theta_Q(X,V_u)
\langle Q,X\rangle-\mathbf i
\langle J_{V_u}(Q),X\rangle
$,
which, after multiplications, result the decomposition 
\begin{equation}
\varphi =
\varphi_{evn_J} +
\varphi_{odd_J} =
\sum_{Q_i}\varphi_{Q_ievn_J}
+\sum_{Q_i}\varphi_{Q_iodd_J},
\end{equation}
where the first function involves all terms having 
even number of $J_{V_u}$ while
for the other one this number is odd. 
For a fixed boundary point $Z$, these
functions depend just on $X$. One can prove, by formula (\ref{circ}),
that these two functions are in two completely independent subspaces
of functions. The proof will be based on the fact that, on a circle $C$,
the first function appears as a trigonometric rational function 
depending on 
$\sin^2 
\beta_{C^\prime}   
\cos^2 \alpha_C
$
while the other function depends irrationally on these terms. This 
independence implies then that both
$
\varphi_{evn_J} (X,Z)
$
and
$
\varphi_{odd_J} (X,Z)
$
must satisfy the Dirichlet condition. In the next step one decomposes
these functions in the form
\begin{eqnarray}
\varphi_{evn_J} (X,Z)=
\varphi_{evn_Jevn_{J^{(b)}}} (X,Z)+
\varphi_{evn_Jodd_{J^{(b)}}} (X,Z),\\
\varphi_{odd_J} (X,Z)=
\varphi_{odd_Jevn_{J^{(b)}}} (X,Z)+
\varphi_{odd_Jodd_{J^{(b)}}} (X,Z),
\end{eqnarray}
where the first term of each function involves all terms having 
even number of $J^{(b)}_{V_u}$, while, for the second one,
this number is odd. Now using the double pole version of (\ref{circ}),
one can see that all 4 functions in this final decomposition fall in 
completely independent subspaces. But, then, all these 4
functions must satisfy the Dirichlet condition. Now we observe that
\begin{eqnarray}
\kappa_{\mathbf Q_F}:\quad
\varphi_{evn_Jevn_{
J^{(b)}}} (X,Z)\quad\to\quad
\varphi_{evn_Jevn_{J^{(b)}}} (X,Z),\\
\kappa_{\mathbf Q_F}:\quad
\varphi_{evn_Jodd_{J^{(b)}}} (X,Z)\quad\to\quad
-\varphi_{evn_Jodd_{J^{(b)}}} (X,Z),\\
\kappa_{\mathbf Q_F}:\quad
\varphi_{odd_Jevn_{J^{(b)}}} (X,Z)\quad\to\quad
\varphi_{odd_Jevn_{J^{(b)}}} (X,Z),\\
\kappa_{\mathbf Q_F}:\quad
\varphi_{odd_Jodd_{J^{(b)}}} (X,Z)\quad\to\quad
-\varphi_{odd_Jodd_{J^{(b)}}} (X,Z),
\end{eqnarray}
which relations are due to 
$
J^{(a)\prime}=
J^{(a)},
J^{(b)\prime}=-
J^{(b)}
$.
Thus 
$\varphi$
and
$\varphi^\prime$
satisfies the Dirichlet condition simultaneously.

For completing this proof only the above mentioned 
{\it Independence Theorems}
should be established. Suppose the contrary, there is a function $\varphi$
which can be represented as linear combinations both of 
$evn_J$- and $odd_J$-type functions: 
\begin{equation}
\label{indep}
\sum_{\tilde Q_r}\varphi_{\tilde Q_revn_J}=\varphi =
\sum_{Q_i}\varphi_{Q_iodd_J}.
\end{equation}
We may suppose that the trigonometric polynomials 
$\cos^{p_j+q_j-s}\sin^s$ are of the same order 
$n=p_j+q_j$ and there are only finite linear combinations in this 
expression. Consider a circle $C$ on which each term in the linear 
combinations appears in the form (\ref{circ}) with the corresponding
constant $\sin\beta_{C_j^\prime}$ and origin (north pole) $O_{Cj}$ on $C$.
(The origin (pole), $O_C$, on $C$ 
is constructed above Theorem \ref{restric}.) 

Next we work just on the right side of (\ref{indep}), involving only the
$odd_J$-type functions. They can be sorted out into classes according
their parameters
$\sin\beta_{C_j^\prime}$ and $O_{Cj}$.
Re-numerate these classes with $m$ in the form
$\{(\beta_{C^\prime}m,O_{Cq})|m=1,2,\dots ,d\}$ 
such that at least one of the
parameters is different for two distinct $m$'s. From the partial sum
determined by a class factor out
$(1-\sin^2 
\beta_{C^\prime m
}   
\cos^2 \alpha_{Cm})^{1\over 2}
$.
Thus, it appears in the one-term-form
$\tilde S_m(1-\sin^2 
\beta_{C^\prime m
}   
\cos^2 \alpha_{Cm})^{1\over 2}
$,
where $\tilde S_m$ is a rational trigonometric function. Furthermore,
\begin{equation}
\varphi 
=\sum_{m=1}^d\tilde S_m
\sqrt{1-\sin^2 \beta_{C^\prime m}\cos^2 \alpha_{Cm}}
=\sum_{m,r=1}^{d,\infty}A_r\tilde S_m
\sin^{2r}\beta_{C^\prime m}\cos^{2r} \alpha_{Cm},
\end{equation}
where
$
\sum_{r=1}^\infty A_rx^r
$
is the Taylor expansion of $\sqrt{1-x}$ about $x=0$. All the coefficients
$A_r$ are non-zeros. Since the left side of (\ref{indep}) is a 
trigonometric rational functions, all terms 
\begin{equation}
\sum_{m=1}^d\tilde S_m
\sin^{2r}\beta_{C^\prime m}\cos^{2r} \alpha_{Cm},
\end{equation}
considered for a fixed $r$, vanish for big numbers $r$, say, if $r> N$.

Now observe that there exist an everywhere dense open subset, $B\subset C$,
such that, for any fixed point $P\in B$ the function values
$
\sin\beta_{C^\prime m}\cos\alpha_{Cm}(P),
$
considered for all $m=1,\dots ,d$
are distinct and, therefore, the $d\times d$-matrix 
$(T_{mk}=
\sin^{2(k+N)}\beta_{C^\prime m}\cos^{2(k+N)}\alpha_{Cm}(P),
$
where $k=1,\dots ,d$, is non-degenerated. Thus equations
$\sum_m T_{mk}
\tilde S_m(P)=0, 
\forall k,$ implies 
$
\tilde S_m(P)=0 
$
on the whole circle $C$. This argument can be repeated in those cases
when  
$
\tilde S_m
$
is substituted by
$
\tilde S_m
\sin^{2r}\beta_{C^\prime m}\cos^{2r}\alpha_{Cm}, \forall r\geq 0
$. 
Therefore, $\varphi =0$, which proves the independence of the
considered function spaces completely.

This argument repeated for 2-pole functions proves the 
desired independence, first, on the subsets
$(X_F^{(a)},X^{(b)})$.
But this statement obviously implies the independence on the whole X-space.

Both independence theorems can be formulated in terms of polynomials
such that, after considering the decompositions, the 
$evn_J$- resp. $odd_J$-type subspaces 
are spanned by functions involving 
even resp. odd number of $J$'s. This observation leads to a simple
establishment of
the independence and Dirichlet-intertwining theorems regarding
changing basis cases as follows. First consider the complex matrix
$c_{ij}(V_u)$
transforming the fixed basis
$
\mathbf Q_F
$
to the changing one,
$
\mathbf Q(V_u).
$
It is obvious that, 
both in the 1-pole and
2-pole cases, 
both type of subspaces 
regarding
the two systems are transformed to each other by the non-degenerated map 
\begin{equation}
\omega :\quad
\mathcal F_{\mathbf Q_F\{p_i,q_i\}}
\to
\mathcal F_{\mathbf Q(V_u)\{p_i,q_i\}}
\quad ,\quad 
\mathcal{HF}_{\mathbf Q_F\{p_i,q_i\}}
\to
\mathcal{HF}_{\mathbf Q(V_u)\{p_i,q_i\}}
\end{equation}
induced by the basis-transformation. This proves both independence
theorems for the changing basis case immediately. Thus we have:
\begin{theorem}
The $\kappa_{\mathbf Q}$ intertwines the Dirichlet condition both in
the fixed, 
$
\mathbf Q_F
$,
and the changing basis,
$
\mathbf Q(V_u)
$,
cases.

The proof is based on the {\bf Independence Theorem} stating 
that the total space
\[
\mathcal{F}_{\mathbf Q_F,n}=
\sum_{\{n=\sum (p_i+q_i)\}}
\mathcal{F}_{\mathbf Q_F
\{p_i,q_i\}},
\]
defined for a fixed $n$ and $Z$, 
is the direct sum of the independent subspaces
$
\mathcal{F}_{\mathbf Q_F,n,evn_J}
$
and
$
\mathcal{F}_{\mathbf Q_F,n,odd_J}
$,
where, after implementing the above described natural decomposition,
the functions from the first resp. second space contain even resp. odd
number of $J$'s. Both of these subspaces further decompose into the
independent subspaces
$
\mathcal{F}_{\mathbf Q_F,n,par_J,evn_{J^{(b)}}}
$
and
$
\mathcal{F}_{\mathbf Q_F,n,par_J,odd_
{J^{(b)}}}
$
defined by the options $par=evn$ or $odd$ given for the parities of 
the number of
$
J^{(b)}
$'s in the expressions. Vector $Z$ should not be the same for the 
participating functions but it can be chosen individually and 
independently both for the even- and odd-type functions.  

The independence guaranties that all 4 component
functions of a $\varphi$ and $\varphi$ itself satisfy the Dirichlet 
condition always simultaneously. The same statements hold for the function
spaces 
$\mathcal{HF}_{\mathbf Q_F,n}$
as well as for both versions of function spaces defined by changing basis
fields.
\end{theorem}

\section{Intertwining of the Neumann conditions.}\label{nc}

The Neumann conditions create a new more complicated situation which
requires reformulations of the proofs given for the Dirichlet conditions
at several points. 
The Z-Neumann conditions, however, which require the vanishing of the 
derivatives of functions taken from the Z-radial directions at the boundary
points, can be strait-forwardly traced back to the Dirichlet conditions.
Let it also be mentioned that the proofs on the boundary manifolds 
exploit only the intertwining of the Dirichlet 
and Z-Neumann conditions. By this reason, the Z-Neumann conditions are 
considered first.
 
In the following computations the integral defining the Z-Fourier
transform is considered on the polar coordinate system. 
The computations are carried out by formulas
$
\partial_{|Z|}
e^{\mathbf i\langle Z,V\rangle}
=
\mathbf i\langle Z_u,V_u\rangle |V|
e^{\mathbf i\langle Z,V\rangle}
=
|Z|^{-1}|V|
\partial_{|V|}
e^{\mathbf i\langle Z,V\rangle}
$
combined with integration by parts. Without loosing the generality, 
one can suppose that the 
test-function, $\phi$, vanishes at the infinity. Then, in terms of
$\phi^\prime :=\partial_{|Z|}\phi$, we have
\begin{eqnarray}
\partial_{|Z|}
\mathcal F_{\mathbf Q\{p_i,q_i\}}(\phi )
=-|Z|
(\mathcal F_{\mathbf Q\{p_i,q_i\}}(
|V|\phi^{\prime} )+
l\mathcal F_{\mathbf Q\{p_i,q_i\}}(\phi )).
\end{eqnarray}
Therefore, a function
$
\mathcal F_{\mathbf Q\{p_i,q_i\}}(\phi )
$
satisfies the Z-Neumann condition if an only if 
$
\mathcal F_{\mathbf Q\{p_i,q_i\}}(
|V|\phi^{\prime} )+
l\mathcal F_{\mathbf Q\{p_i,q_i\}}(\phi )
$
satisfies the Dirichlet condition. Since the Dirichlet condition
is intertwined in all cases, we have
\begin{theorem}
The Z-Neumann condition is intertwined both in the fixed and changing
basis cases.
\end{theorem}

From now on, the standard Neumann condition is scrutinized.
By formulas (3.2) and (3.7)
of \cite{sz2}, the normal vector at a boundary point $(X,Z)$ 
and the Laplacian on the boundary manifolds are of the
form 
\begin{eqnarray}
\mathbf\mu =A(|X|,|Z|)X_u+B(|X|,|Z|)J_Z(X)+C(|X|,|Z|)Z_u,  
\label{normv}\\
\tilde\Delta=\Delta_{S_X(Z)}+(1+\frac 1{4}|X|^2)\Delta_{S_Z(X)}
+\sum_{\alpha =1}^{l-1}\partial_\alpha D_\alpha \bullet,
\label{bound_lapl}
\end{eqnarray}
where $S_X(Z)$ is the X-sphere over $Z$ and $S_Z(X)$ is the Z-sphere
over $X$, furthermore, $\{\partial_1,\dots ,\partial_{l-1}\}$ 
is an orthonormal basis in the tangent space of the Z-sphere $S_Z(X)$
at $Z$. On sphere$\times$ball-type manifolds 
the first term of (\ref{normv}) should be omitted. 

The standard Neumann condition is considered, first, for  
one pole functions. 
The derivative with respect to the normal direction
is built up by X- resp. Z-radial-derivatives and $J_Z(X)\bullet$. 
If the foot of perpendicular through $X$ to $S_Q$ is $X_Q$, and thus 
$X=X_Q+X_Q^\perp$ holds, then 
$J_Z(X)\bullet =J_Z(X_Q)\bullet +J_Z(X^\perp_Q)\bullet$. 
First  
$
J_{Z}(X_Q)\bullet  
$
is considered.
If $X_Q=|X_Q|(\cos (\beta )Q+\sin (\beta )J_{Z_0}(Q) \in P(Q,Z_0)$, then
\begin{equation}\label{nc1}
J_Z(X_Q)=|X_Q|(\cos(\beta )J_Z(Q)-
\sin (\beta ) 
\langle Z_0,Z\rangle Q
+\sin (\beta )J_{Z^\perp}J_{Z_0}(Q)), 
\end{equation}
where $Z^\perp$ is the perpendicular component of $Z$ to $Z_0$. 
Therefore, the 
$J_{Z^\perp}$ and $J_{Z_0}$ 
are anti-commuting and the
$J_{Z^\perp}J_{Z_0}$ is a skew endomorphism. Thus,
\begin{eqnarray}\label{nc2}
\langle Q,J_Z(X_Q)\rangle=-|X_Q|\sin\beta \langle Z_0,Z\rangle\, ,\,
\\
\langle J_{Z_0}(Q),J_Z(X_Q)\rangle=|X_Q|\cos\beta 
\langle Z_0,Z\rangle ,
\\
J_Z(X_Q)\bullet
\cos\alpha =-|X_Q|\sin\beta \langle Z_0,Z\rangle\, ,\,
\\
J_Z(X_Q)\bullet\sin\alpha=|X_Q|\cos\beta \langle Z_0,Z\rangle ,
\\
J_Z(X_Q)\bullet
\mathcal {HF}_{Qpq}(\phi )(X,Z)=
\sum
|X_Q|^{p+q}
\cos^{p+q-s}(\beta )\sin^s(\beta )
\tilde S^T_{spq},
\\
{\text where}\,\,
\tilde S^T_{spq}=-|X_Q|
\langle Z_0,Z\rangle
((p+q-s+1)\tilde P_{(s-1)pq}-
(s+1)\tilde P_{(s+1)pq}).  
\label{nc7}
\end{eqnarray}

The computations with $J_Z(X_Q^\perp)\bullet$ are based on
\begin{eqnarray}\label{nc8}
\langle Q,J_Z(X^\perp_Q)\rangle=0\quad ,\quad
\langle J_{Z_0}(Q),J_Z(X^\perp_Q)\rangle=
-\langle J_{Z^\perp}J_{Z_0}(Q),X_Q^\perp\rangle .
\end{eqnarray}
Since differentiation $J_Z(X^\perp_Q)\bullet$ acts, non-trivially, 
on the considered 
functions only by its contribution
$
\langle J_{Z_0}(Q),J_Z(X^\perp_Q)\rangle(\cos\beta\partial_\alpha
+\sin\alpha\partial_r) 
$
to the $\partial_\alpha$- and the radial $\partial_r$-direction,
therefore:
\begin{eqnarray}
J_Z(X^\perp_Q)\bullet\cos\alpha =
\langle J_{Z^\perp}J_{Z_0}(Q),X_Q^\perp\rangle\cos\beta\sin\beta ,\\
J_Z(X^\perp_Q)\bullet\sin\alpha=
-\langle J_{Z^\perp}J_{Z_0}(Q),X_Q^\perp\rangle
\cos\beta\cos\beta ,
\\
J_Z(X^\perp_Q)\bullet
\mathcal {HF}_{Qpq}(\phi )(X,Z)=
|X_Q|^{p+q}\sum_{s}
\cos^{p+q-s}(\beta )\sin^s(\beta )\tilde S^\perp_{spq} ,\\
\tilde S^\perp_{spq}=
-\langle J_{Z^\perp}J_{Z_0}(Q),X_Q^\perp\rangle
((p+q)|X_Q|^{-1}\sin\beta\tilde P_{spq}-
\cos\beta \tilde D_{spq}),\\
\tilde D_{spq}=(p+q-s+1)\tilde P_{(s-1)pq}-
(s+1)\tilde P_{(s+1)pq}.  
\end{eqnarray}

A preliminary version of the standard Neumann condition 
for a one-pole function can be stated in the
following form. A $\phi$-generated one-pole function satisfies the 
Neumann condition if and only if
\begin{eqnarray}\label{nc14}
\tilde R_{spq}=(A\partial_{|X|}
+C\partial_{|Z|})\tilde P_{spq}
+B(\tilde S^T_{spq}+\tilde S^\perp_{spq})=0 
\end{eqnarray}
holds at the boundary points, for all $Z_0(X)$ and $0\leq s\leq p+q$, 
where coefficients $A,B$ and $C$ are defined by the normal 
vector $\mu$. Already this version reveals that only the intertwining
regarding the first two terms in this condition can be traced back to
the intertwining of the Dirichlet conditions. The other two terms
involve functions such as
$
\langle {Z_0},Z\rangle
$
and
$
\langle J_{Z_0}(Q),J_Z(X^\perp_Q)\rangle=
|X^\perp_Q|
\langle J_{Z_0}(Q),J_Z(
X^\perp_{Q0}
)\rangle
$
which appear outside of the integral terms. It is noteworthy that the
second function is zero on spaces 
$H^{(a,b)}_3$
and 
$H^{(a,b)}_7$. This is due to the fact that
the irreducible components 
$H^{(1,0)}_3$
and 
$H^{(1,0)}_7,$
yield the well known $J^2$-condition, meaning, that 
for any product 
$J_{Z^\perp}J_{Z_0}$
there exist $J_{\tilde Z}$
such that
$J_{Z^\perp}J_{Z_0}=J_{\tilde Z}$ holds. 
On arbitrary H-type groups, for fixed $Z$ and 
unit vector
$
X^\perp_{Q0},
$
all these terms define polynomials which are suitable to establish
the independence theorems seen for the Dirichlet condition also for
the Neumann condition.
Similar formulas can be established also for the two-pole functions
which also yield the corresponding independence theorem. Finally we get:
\begin{theorem}
The $\kappa_{\mathbf Q}$ intertwines the standard Neumann conditions 
both in the fixed, 
$
\mathbf Q_F
$,
and the changing basis,
$
\mathbf Q(V_u)
$,
cases.

The proof is based on observing that the total space
\[
\mu\bullet \mathcal{HF}_{\mathbf Q_F,n}=\mu\bullet
\sum_{\{p_i,q_i\}}
\mathcal{HF}_{\mathbf Q_F
\{p_i,q_i\}},
\]
defined for fixed $Z$, $n$ and running  
$\{p_i,q_i\}$ satisfying
$
n=\sum (p_i+q_i)
$,
is a direct sum of independent subspaces
$
\mathcal{HF}_{\mathbf Q_F,n,evn_J}
$
and
$
\mathcal{HF}_{\mathbf Q_F,n,odd_J}
$,
which are defined such that
the functions from the first resp. second space contain even resp. odd
number of $J$'s (i. e. , $\sin\beta$'s, according to the above 
decomposition). Both subspaces further decompose into the
independent subspaces
$
\mu\bullet\mathcal{F}_{\mathbf Q_F,n,par_J,evn_{J^{(b)}}}
$
and
$
\mu\bullet\mathcal{F}_{\mathbf Q_F,n,par_J,odd_
{J^{(b)}}}
$
defined by the options $par=evn$ or $par=odd$, available 
for the parity of number of
$
J^{(b)}
$'s in the expressions. The independence guaranties that all the 4 component
functions of a $\varphi$ along with $\varphi$ satisfy the 
standard Neumann 
condition always simultaneously. The same statements hold for the function
spaces 
$\mathcal{F}_{\mathbf Q_F,n}$
as well as for both versions of function spaces defined by changing basis
fields.
\end{theorem}

\section{Intertwining on the boundary manifolds. 
\label{bound_int}}

Since the intertwining operator preserves the Dirichlet condition,
by restrictions, it induces a well defined bijection between the 
$L^2$ spaces defined on the boundaries. 
Since each smooth function on the boundary
extends to ones satisfying 
the Z-Neumann condition, furthermore,
this condition is also preserved by the operator,  
it is enough to represent the functions on the boundary  by
restrictions of those satisfying the Z-Neumann condition.  

If $\partial_l=\partial_{|Z|}$ is the Z-partial derivative with respect to
the normal direction $Z_u$, then the angular momentum operator
$\mathbf M$ (resp. $\tilde{\mathbf M}$) 
on the ambient (resp. boundary) manifold differ from each other just by
$\partial_{|Z|}D_{Z_u}\bullet$. This operator vanishes on functions
satisfying the Z-Neumann condition, thus  
$\mathbf M$ and $\tilde{\mathbf M}$ acting on these functions provide
the same results. This
argument proofs that not just 
$\mathbf M$ and $\mathbf M^\prime$ 
but also
$\tilde{\mathbf M}$ and $\tilde{\mathbf M}^\prime$ 
are intertwined by the operator. This is the most crucial part in the
proof of the pursued theorem. The intertwining regarding 
$\Delta_{S_X(Z)}$ 
has already been established on the ambient manifold, thus one should
consider only  
$\Delta_{S_Z(X)}$. Since the intertwining regarding $\Delta_Z$ is
established on the ambient manifold, only the intertwining of
the radial Laplacian $\Delta_{|Z|}$ should be established. Since
this operator acts on functions satisfying the Z-Neumann condition,
the question is if 
$\partial^2_{|Z|^2}$ 
is invariant under
the action of the operator. This statement immediately follows from
the following computations where the integral defining the Z-Fourier
transform is considered on the polar coordinate system. 
The computations start out with
$
\partial^2_{|Z|^2}
e^{\mathbf i\langle Z,V\rangle}
=
-\langle Z_u,V_u\rangle^2|V|^2
e^{\mathbf i\langle Z,V\rangle}
=
|Z|^{-2}|V|^2
\partial^2_{|V|^2}
e^{\mathbf i\langle Z,V\rangle}
$
and are completed by integration by parts. Then, in terms of
$\phi^\prime =\partial_{|Z|}\phi$, we have
\begin{eqnarray}
\partial^2_{|Z|^2}
\mathcal F_{Qpq}(\phi )
=|Z|^{-2}(
\mathcal F_{Qpq}(
|V|^2\phi^{\prime\prime} )+
2(l+1)\mathcal F_{Qpq}(
|V|\phi^\prime )+
l\mathcal F_{Qpq}(\phi )),
\end{eqnarray}
where the $|Z|=R_Z$ is a constant. 
Like in case of the Z-Neumann condition, this formula establishes the
desired intertwining property for
$\partial^2_{|Z|^2}$.

\section{Intertwining on solvable extensions.}

The isospectrality theorems naturally extend to 
the solvable extensions. 
The Laplacians on the ambient- and boundary-manifolds,
furthermore, the normal vectors at the boundaries are described
in formulas (1.12), (3.30), and (3.29) of \cite{sz2}.
The generator functions
are of the form $\phi(|X|,t,V)$ in this case, till, the intertwining 
operator is defined by the same Z-Fourier transform like on H-type groups.
I. e., the 
intertwining on the solvable group is completely determined by its action
induced on the nilpotent group. The details in \cite{sz2} show that the
terms due to the t-variable, which are in a certain combination with the
terms of the Laplacian defined on the nilpotent group, make no effect
on proving the intertwining for the solvable groups in the same way as 
for the H-type groups.

\section{The new striking examples.}

There is a subgroup, 
$\mathbf {Sp}(a)\times \mathbf{Sp}(b)$, 
of isometries
on a Heisenberg-type group $H^{(a,b)}_3$ which acts as the identity
on the Z-space.
Note that 
these isometries act transitively on the
X-spheres 
of $H_3^{(a+b,0)}$. 

The complete isotropy group of isometries fixing the origin is
$(\mathbf {Sp}(a)\times \mathbf{Sp}(b))\cdot SO(3)$,
where the action of $SO(3)$, described in terms of unit quaternions
$q$ by 
$
\alpha_q(X_1,\dots ,X_{a+b},Z)=
(qX_1\overline q,\dots ,qX_{a+b}\overline q,qZ\overline q),
$ 
is transitive on the Z-sphere. The elements of this isotropy group induce 
isometries on the 
sphere$\times$sphere-type submanifolds, furthermore, there is proved in the
Extension Theorem of \cite{sz2}
that these are the only isometries on these submanifolds. Note that these
isospectral manifolds in a family have non-isomorphic isometry groups of
different dimensions such that
they are homogeneous in
$H_3^{(a+b,0)}\cong H_3^{(0,a+b)}$, 
while the other members are locally inhomogeneous.

The sphere$\times$sphere-type submanifolds of the solvable extensions of
H-type groups is defined such that, 
over each point of a sphere in the X-space, 
one considers the same geodesic sphere around the origin $(0,1)$ 
of the hyperbolic $(Z,t)$-space. 
The isotropy group of isometries acting on 
$SH^{(a,b)}_3$, where $ab\not =0$, is 
$(\mathbf {Sp}(a)\times \mathbf{Sp}(b))\cdot SO(3)$,
while 
it is $\mathbf{Sp}(a+b)\cdot\mathbf{Sp}(1)$ 
on 
$SH^{(a+b,0)}_3\cong SH^{(0,a+b)}_3$. 
The above statements extend also to these groups. 

These are
the new striking examples brought by the reconstructed intertwining
operator. The sphere-type striking examples discussed in the earlier
papers can be explained similarly. They are the geodesic spheres defined
by the same radius for the family
$SH^{(a,b)}_3$. In this case the metric on
$SH^{(a+b,0)}_3$ is two-point homogeneous, therefore, having homogeneous
geodesic spheres. The geodesic spheres on the other groups are locally
inhomogeneous. 

It should be mentioned that Sch\"uth \cite{sch1} constructed the 
literature's first
isospectral metrics defined on simply connected manifolds on Cartesian
products of spheres. Among them there are also sphere$\times$sphere-type
manifolds. Her construction arises from a completely different setting,
however, and all provided metrics are locally inhomogeneous.

\section{Isospectralities for $\sigma$-equivalent metrics.}
 
A $\sigma$-deformation of an
endomorphism space 
$J_{\mathbf z}$ 
is defined by an 
involutive orthogonal transformation,
$\sigma$, of the X-space which commutes with all endomorphisms
of $J_{\mathbf z}$. The $\sigma$-deformed endomorphism space consists
of endomorphisms $\sigma J_Z$. By using 
irreducible decomposition regarding 
the orthogonal Lie algebra generated by the elements of 
$
J_{\mathbf z}$, for any $\sigma$-deformation, there exist a decomposition
$
\mathbf v=\mathbf v
^{(a)}
\oplus\mathbf v^{(b)}
$
with components invariant under the action of the endomorphisms such
that for their restriction onto the components the relations
$
\sigma J_Z^{(a)}=J_Z^{(a)}
$
and
$
\sigma J_Z^{(b)}=-J_Z^{(b)}
$
hold. Note that the family, 
$J^{(a,b)}_l$, 
of Cliffordian endomorphism spaces 
defined by the same $a+b$ and $l$
consists of $\sigma$-equivalent endomorphism spaces. In papers
\cite{sz1,sz2} the isospectrality is stated on the ball- and 
sphere-type domains of such $\sigma$-equivalent 2-step nilpotent Lie 
groups and their 
solvable extensions whose endomorphism spaces 
contain at least one anticommutator. Thus the extension of the
isospectrality theorems to
$\sigma$-equivalent 2-step nilpotent Lie groups and their solvable
extensions provide plenty additional examples to those produced by the
anticommutator technique. 
For the sake of simplicity, we consider such
two step nilpotent Lie groups whose endomorphism spaces contain at least
one non-degenerated endomorphism. Then almost all endomorphisms must be 
non-degenerated acting on an even dimensional X-space.
 
First, we look for the necessary modifications which make the techniques 
developed for H-type groups working also for 
$\sigma$-deformations.
The Laplacians
in these general cases differ from the Laplacians of
H-type groups just by the term
$(1/4)|X|^2\Delta_Z$, what is now
$(1/4)\sum\langle J_{\alpha}(X),J_{\beta}(X)\rangle
\partial^2_{\alpha\beta}$. 
Also the intertwining operator, defined for the changing basis case, 
must be modified as follows. Let 
$(Q_1(V_u),\dots ,Q_{k/2}(V_u))=(\mathbf Q^{(a)},\mathbf Q^{(b)})$ 
be an appropriate orthonormal basis 
field such that each $Q_{V_ui}$ is an eigenvector of 
$J^2_{V_u}$ with eigenvalue $-\lambda_i^2(V_u)$.
Let $\tilde J_{V_u}$ be the normalized endomorphism which has the same
kernel as $J_{V_u}$ and is defined by $(1/\lambda_i(V_u))J_{V_u}$
on the maximal eigensubspaces belonging to $\lambda_i>0$. Note that 
the kernel is trivial for an everywhere dense open subset of the unit 
vectors $V_u$, furthermore, 
this endomorphism may not be in $J_{\mathbf z}$. 
Then, by definition,
$\Theta_{Q(V_u)}(X,V_u)=\langle Q +\mathbf i\tilde J_{V_u}(Q),X\rangle$. 
The complex coordinate system defined by this basis for non-degenerated
endomorphisms is denoted
by $\{z_{V_u1},\dots ,z_{V_uk/2}\}$. Then the intertwining 
is defined by these changing complex 
coordinates such that functions $\phi$ 
should be of the form $\phi (|X_{V1}|,\dots ,|X_{Vr}|,V)$, 
where $X=\sum X_{Vi}$ is the decomposition regarding the 
eigenspaces of $J_V^2$. 
This requirement is slightly different from that what is 
considered on H-type 
groups, but one can reach to them in the same way: Start with
functions depending just on $V$, first. Then, it turns out 
that the same operator is defined by the above more complicated functions.
Furthermore, the domain is the largest possible, containing the complete
$L^2$- function space. Also Theorem \ref{euclapl} concerning the 
intertwining of the Euclidean Laplacian and radial functions 
on the X-space, remains true for $\sigma$-deformations. Therefore,
by $J^{\prime 2}_V=\sigma J_V\sigma J_V=\sigma^2J_V^2=J_V^2$ 
and formulas (\ref{MF}), (\ref{Delta_ZF}), 
which also extend to $\sigma$-deformations,
this is indeed an operator intertwining the Laplacians
$\Delta$ and $\Delta^\prime$ term by term.

One should check out also the intertwining of the boundary conditions.
First note that the technique developed for H-type groups works out
strait-forwardly only on groups where the intertwining can be established 
also  with a fixed basis 
$\mathbf Q_F$,  
therefore, one can use one- and two-pole functions with poles being in
$span_{\mathbb R}(\mathbf Q_F)$. 
Also well-defined polar coordinate systems established on 
$\mathbb{R}Q\oplus\tilde J_{\mathbf z}(Q)$, 
where 
$\tilde J_{\mathbf z}$
is spanned by endomorphisms of the form
$\tilde J_{V_u}$,
are necessary conditions for this technique.
All these requirements are satisfied only in those cases where, for all
unit pole $Q\in span_{\mathbb R}(\mathbf Q_F)$, the unit vectors 
$\tilde J_{\mathbf z}(Q)$ 
form an everywhere dense open subset of the unit sphere of the space
spanned by these vectors (equator). Then, any such vector is connected 
with the pole by an
$\alpha$-parameter circle defined for $0\leq\alpha\leq\pi$ such that it
has the parameter $\pi /2$. The possible missing circles, 
which are due to the 
degenerated endomorphisms, can be implemented by 
limiting. The extension to this cases works out after
other additional modifications. 
   
First we check on formula (\ref{theta}) of Section \ref{tdc}, where,
on a parameter-circle, the corresponding $Z_0$ should be exchanged for
$Z_*(Z_0)$ defined by the dual
of the functional $\varphi (Z)=\langle\tilde J_{Z_0}(Q),J_Z(Q)\rangle$. 
Yet, the $Z_*(Z_0)$ is a polynomial function of $Z_0$,
implying that functions $\tilde R_{spq}$, introduced in (\ref{circ}),
will be polynomials on the above equator. 
Thus the computations can be processed in the same way as  
earlier. Since for $\sigma$-deformations the relations 
$\varphi (Z)=\varphi^\prime (Z)\, , Z_*=Z^\prime_*$ hold,
the proof regarding the Dirichlet and Z-Neumann condition can be 
completed by the same argument seen for H-type groups.
 
Regarding the Neumann condition the $Z_0$ inside of the integral term
should be exchanged for $Z_*$, while terms 
$\langle Z_0,Z\rangle$
resp. 
$\langle J_{Z^\perp}\tilde J_{Z_0}(Q),X^{\perp}_Q\rangle$
outside of the integral
should be exchanged for much more complicated expressions. Even so,
they provide polynomial functions and
the modified computations concerning
formulas (\ref{nc2})-(\ref{nc7}) along with  
$
\langle J_{Z^*}\tilde J_{Z_0}(Q),Q^*\rangle =
\langle J^\prime_{Z^*}\tilde J^\prime_{Z_0}(Q),Q^*\rangle 
$ 
yield the intertwining also of the Neumann conditions for the 
$\sigma$-deformations whose endomorphism spaces satisfy the above 
conditions. In these cases, the theorem extends also to 
the boundary manifolds and the solvable extensions.

Fortunately, one should not go through the steps of this complicated 
construction which is incomplete without scrutinizing the question
of the existence. 
Basically, what the above proof exploits is that 
the operator introduced by the
changing-basis-technique for $\sigma$-deformations intertwines 
both the Laplacians and
boundary conditions if the manifold satisfies the independence theorems.
But this theorems are certainly yielded on groups having endomorphism
spaces which are the results of ``small'' perturbations performed on the
endomorphism space of H-type groups such that, by choosing a new
endomorphism space close to the Clifordians, one perturbs the 
endomorphism space acting on the irreducible space $\mathbb R^{r(l)}$.
This defines uniquely determined 
perturbations for the reducible endomorphisms (see more details in the
end of Introduction). If the perturbation is ``small''
which changes the endomorphisms just slightly, the subspaces 
in the independence theorems 
(which are closed in the ambient Hilbert space) 
keep being independent. Even the conditions for the existence
of polar coordinate systems
are satisfied on manifolds defined by ``smaller'' perturbations.
However, the above construction can be left out completely 
because the idea of perturbation 
provides more examples than those provided by the above process. 
Also the non-isometry 
proofs, guaranteeing that the considered isospectral metrics have 
different
local geometries, are inherited from the metrics being ``slightly'' 
perturbed. For the latter metrics the non-isometry proofs are completely
established in \cite{sz1,sz2}.  
Thus we have:
\begin{theorem} 
The intertwining operator introduced 
for $\sigma$-equivalent groups   
by the changing-basis-technique 
always intertwines the Laplacians both on ball$\times$ball-
and sphere$\times$ball-type manifolds. Furthermore,
there exist an open neighborhood, $U$, in the space of $l$-dimensional
endomorphisms spaces with endomorphisms acting on the irreducible space 
$\mathbb R^{r(l)}$ of a given Clifford
endomorphism space, $J_{\mathbf z}$ such that the latter one 
is in $U$ and all 2-step nilpotent groups and
their solvable extensions constructed
by endomorphism spaces belonging to $U$ satisfy the independence theorems.
Therefore, for these groups, the intertwining operator intertwines also
the boundary conditions. In other words, small perturbations performed 
on the endomorphism spaces of H-type groups provide a wide range of 
$\sigma$-equivalent groups which are isospectral on the 
corresponding ball$\times$ball- and sphere$\times$ball-type submanifolds.  

This theorem extends to solvable groups as well as to the boundary 
manifolds, both in the nilpotent and solvable cases. For dimensions 
satisfying $l=4n+3$ also such $U$ exist for which the corresponding
metrics in an isospectrality family have different local geometries.
\end{theorem} 

\noindent{\bf Remark.} The key point in the above process dealing with 
one-pole functions is that the multi-linear function 
\begin{equation}
H(X,X^*,Z,Z^*):=\langle J_Z(X),J_{Z^*}(X^*)\rangle
\end{equation} 
does not change
during $\sigma$-deformations. 
The Ricci tensor can be described in terms of this function by
\begin{eqnarray}
R(X,X^*)=-(1/2)\sum H(X,X^*,Z_\alpha ,Z_\alpha)\, ,\,
\\ 
R(Z,Z^*)=(1/4)\sum H(E_i,E_i,Z,Z^*) \quad ,\quad R(X,Z)=0
\end{eqnarray}
(cf. formula (1.9) of \cite{sz1}),
thus also this tensor is not changing during $\sigma$-deformations.

The Gordon-Wilson \cite{gw} isospectrality examples were constructed
on ball$\times$torus-type manifolds 
by spectrally equivalent endomorphism spaces, meaning the existence
of orthogonal transformations associating isospectral
endomorphisms to each other. More precisely, they constructed 
continuous families of metrics which are isospectral on functions. 
These metrics are not isospectral on 1-forms, 
however, due to the fact that
the norm of the Ricci tensor is changing during these
deformations \cite{sch2}. The question arises if the 
domains investigated in
this paper are isospectral on the Gordon-Wilson examples.

Even though the changing Ricci tensor strongly suggests the negative
answer, the question is more complicated. 
Indeed, the intertwining operator,
constructed on ball$\times$torus-type manifolds such that
the globally defined operator is the direct sum of operators
constructed on the invariant subspaces
$W_\gamma$
by the single endomorphism $J_\gamma$ separately, 
provides an operator intertwining
the Laplacians also for the
ball$\times$ball- and sphere$\times$ball-type manifolds which are in
the ball$\times$torus-type manifold. Actually, this is a discrete version
of those constructions where one is using changing basis. 

There are many problems arising
when one tries to establish the intertwining of the boundary conditions
for this operator on the Gordon-Wilson examples. 
First of all, one can not pass to a fixed basis and involve one- and 
two-pole functions because no fixed basis exists which is transformed
to a well defined fixed basis by all those point transformations which 
define 
the intertwining operator for the subspaces $W_\gamma$. The metrics in the
Gordon-Wilson examples are out of the touch also of the perturbation 
technique. 
Even though the independence theorems were guaranteed in an other way,
they would work together just with the $\sigma$ deformations for which the
decomposition 
$
\mathbf v= 
\mathbf v^{(a)}\oplus 
\mathbf v^{(b)}
$ 
holds. (This argument shows that the complete isospectrality can not be
directly established by this discrete version of the intertwining
operators even in case of $\sigma$-deformations, because,it is necessary
to involve 1-and 2-pole functions defined by a constant basis.)  
In short, there is no way to prove the intertwining of the
boundary conditions by our technique. 

On the other hand, by the above
considerations, a non-changing Ricci tensor is 
always a necessary condition
for intertwining even the Dirichlet conditions
on ball$\times$ball-type manifolds. Therefore, 
these particular continuous
operators change the Dirichlet conditions along with the Dirichlet 
spectra. Thus, they can not induce  
operators transforming functions defined on the boundaries either.
In other words, operators defined by restrictions onto the boundary
manifolds are not well defined regarding the GW-deformations. 
Even though they were introduced
by suitable reductions, they would change the spectra also on the
boundary manifolds. This phenomena strongly suggests that 
no other suitable operators exist and the spectra of the
ball$\times$ball- and sphere$\times$ball-type manifolds change during
these continuous deformations.  

\noindent{\bf Acknowledgements.} 
The author is indebted for the hospitality
and excellent working conditions provided by the Max Planck Institute 
for Mathematics in the Sciences, Leipzig, in the academic year 2007/2008, 
where important key ideas of this paper were found.

\end{document}